%
%

\def\dateofversion{April 22, 2001} 

%
%

\def\picturelist{abd-draw*.eps (* ranging from A to G)}

%
%

\def\formato{1}    
                   %
\def\stampa{0}     
                   %
\def\nicefonts{0}  
                   %
\def\insertfig{2}  
                   %
\def\dialogo{0}    
                   
%
%

\nopagenumbers

\ifnum\formato=1
   \magnification 1100
   \voffset 2 cm
   \hoffset 1 cm
\else
   \magnification 1320
   \voffset 0.2 cm
   \hoffset -0.1 cm
   ~\vfill\eject
\fi

\tolerance 1000
\vsize 18.5 cm
\hsize 12.4 cm
\parindent 15 pt

%
%

\pageno=1

\def\leftheadline{\eightrm\folio\hfill G. Alberti, G. Bouchitt\'e, 
G. Dal Maso \hfill\phantom\folio}

\def\rightheadline{\eightrm\phantom\folio\hfill The calibration 
method for the Mumford-Shah functional\hfill\folio}

\def\nullheadline{\eightrm Revised version, \dateofversion \hfill}

\headline{%
\ifnum\pageno=1\nullheadline\else%
\ifodd\pageno\rightheadline\else%
\leftheadline\fi\fi%
}

%
%
\newlinechar=`@

\ifnum\dialogo=0

\message{
@
@ A T T E N T I O N !!
@
@ If some fonts or pictures files are not available, 
@ try to change the control parameters listed 
@ at the beginning of this file.
@
@
}

\else

\message{Choose one of the following typeset modes:
@  1. COMPLETE / Textures on Mac
@  2. COMPLETE / standard TeX on UNIX/PC/Mac
@  3. MINIMAL
@------------------------------------------------------------
@[1/2] require the fonts rsfs and msbm, and macroes for
@  handling eps files (epsf.def and epsf.sty, resp.);
@  the files of the pictures (\picturelist)
@  should be placed in the same directory as the TeX source.
@[1] warning: pictures may not appear in preview!
@[3] pictures will not be included in the output!
@------------------------------------------------------------
@[type 1/2/3]: }
\read-1 to\answer
\ifnum\answer=3
    \def\nicefonts{0}
    \def\insertfig{0}
\else
    \def\nicefonts{1}
    \ifnum\answer=2\def\insertfig{2}\else\def\insertfig{1}\fi
\fi

\fi

%
%

\def\NOTA#1{\ifnum\stampa=2\medskip\hrule\smallskip{\sl #1}\smallskip
\hrule\bigskip\else\fi}

%
%

\ifnum\insertfig=0
   \def\figureps#1#2{\vbox{
       \vskip .7 cm
       \centerline{FIGURE [#2] TO BE INSERTED HERE}
       \vskip .7 cm
       \centerline{\paragrafo Figure #1}\medskip
     }}
\else
   \def\figureps#1#2{\vbox{
       \centerline{\epsfbox{#2}}
       \centerline{\paragrafo Figure #1}\medskip
     }}
\fi

\ifnum\insertfig=1
   \input epsf.def
\else\ifnum\insertfig=2
   \input epsf.sty
\else\fi\fi

%
%

%

\font\sixrm=cmr6

\newcount\chapno\chapno=0          
\newcount\tagno \tagno=0           
\newcount\thmno \thmno=0           
\newcount\bibno \bibno=0           
\newcount\verno                    

\newif\ifproofmode

\ifnum\stampa=0 \proofmodefalse\else\proofmodetrue\fi

\newif\ifwanted\wantedfalse
\newif\ifindexed\indexedfalse


\def\ifundefined#1{\expandafter\ifx\csname+#1\endcsname\relax}

\def\Wanted#1{\ifundefined{#1} \wantedtrue
\immediate\write0{Wanted #1
\the\chapno.\the\thmno}\fi}

\def\Increase#1{{\global\advance#1 by 1}}


\def\Assign#1#2{\immediate
\write1{\noexpand\expandafter\noexpand\def
    \noexpand\csname+#1\endcsname{#2}}\relax
    \global\expandafter\edef\csname+#1\endcsname{#2}}

\def\lPut#1{\ifproofmode\llap{\hbox{\sixrm #1\ \ \ }}\fi}
\def\rPut#1{\ifproofmode$^{\hbox{\sixrm #1}}$\fi}


\def\chp#1{\global\tagno=0\global\thmno=0\Increase\chapno
\Assign{#1}
{\the\chapno}{\lPut{#1}\the\chapno}}


\def\thm#1{\Increase\thmno
\Assign{#1}{\the\chapno.\the\thmno}\the\chapno.\the\thmno
       \rPut{#1}}


\def\frm#1{\Increase\tagno
     \Assign{#1}{\the\chapno.\the\tagno}
       \lPut{#1}{\the\chapno.\the\tagno}}


\def\bib#1{\Increase\bibno
\Assign{#1}{\the\bibno}\lPut{#1}{\the\bibno}}


\def\rf#1{\Wanted{#1}\csname+#1\endcsname\relax\rPut {#1}}

\input A.sty
\Increase\verno
\immediate\openout1=A.sty
\immediate\write1{\noexpand\verno=\the\verno}

\ifindexed
\immediate\openout2=\jobname.idx
\immediate\openout3=\jobname.sym
\fi

%
%

\newfam\Rsfs
\newfam\Msbm

\font\eightti=cmti8
\font\eightrm=cmr8
\font\eightbf=cmbx8
\font\eighttt=cmtt10 scaled 800

\font\titolo=cmsl12 scaled 1200
\font\autore=cmcsc10
\font\paragrafo=cmcsc10
\font\sezione=cmbx12
\font\sottosezione=cmsl10 scaled 1100

\ifnum\nicefonts=0
   \def\Bf{\bf}
   \def\Cal{\cal}
\else
   \font\trea=msbm10
   \font\treb=msbm7
   \font\trec=msbm5
   \textfont\Msbm=\trea
   \scriptfont\Msbm=\treb
   \scriptscriptfont\Msbm=\trec
   \def\Bf{ \fam\Msbm \trea}

   \font\duea=rsfs10
   \font\dueb=rsfs7
   \font\duec=rsfs5
   \textfont\Rsfs=\duea
   \scriptfont\Rsfs=\dueb
   \scriptscriptfont\Rsfs=\duec
   \def\Cal{\fam\Rsfs \duea}
\fi

%
%

\def\R{{\Bf R}}

\def\SS{{\Bf S}}
\def\L{{\Cal L}}
\def\H{{\Cal H}}
\def\F{{\Cal F}}
\def\G{{\Cal G}}

\def\ll{\lambda}
\def\O{\Omega}
\def\GM{\Gamma}
\def\eps{\varepsilon}
\def\bd{\partial}
\def\loc{_{\rm loc}}
\def\div{{\rm div}}
\def\osc#1{{\rm osc}\,#1}
\def\aplim{\mathop{\rm ap\,lim}}
\def\dps{\displaystyle}
\def\ove{\overline}

\def\nus#1{\mathop{\nu{\trait{0}{0}{3}}}\nolimits_{#1}}


\def\itemm#1#2{{\vskip 1 pt\parindent 37 pt\item{\hbox to 
17 pt{\rm#1\hfill}}#2\par}}

\def\Section#1{\goodbreak\vskip .8 cm\vbox{
\parindent 15 pt\hskip - 15 pt\sezione #1}\vskip .3 cm\nobreak}

\def\subsection#1{\goodbreak\vskip .6 cm\noindent{\sottosezione #1}%
\vskip .3 cm\nobreak}

\def\Theor#1{\medskip{\paragrafo Theorem #1.~--}}
\def\Lemma#1{\medskip{\paragrafo Lemma~#1.~--}}
\def\Def#1{\medskip{\paragrafo Definition~#1.~--}}
\def\Rem#1{\medskip{\paragrafo Remark~#1.~--}}
\def\Pr{\medskip{\paragrafo Proof.~}}
\def\Parag#1{\medskip{\paragrafo #1.~--}}

\def\trait#1#2#3{\vrule width #1pt height #2pt depth #3pt}
%
\def\qed{~~~~\trait{.3}{6}{0}\kern-.3pt\trait{6}{.3}{0}%
\kern-6pt\trait{6}{6}{-5.7}\kern-.3pt\trait{0.3}{6}{0}\goodbreak}
%
\def\vspazio{\vphantom{\trait{0}{13}{0}}}
%
\def\LM#1{\hbox{\vrule width.2pt \vbox to#1pt{\vfill
\hrule width#1pt height.2pt}}}
%
\def\LL{{\mathchoice {\>\LM7\>}{\>\LM7\>}{\,\LM5\,}{\,\LM{3.35}\,}}}
%

%
%

\null
\vskip 1 cm
\centerline{\titolo
The calibration method for the Mumford-Shah functional}
\medskip
\centerline{\titolo
and free-discontinuity problems}

%
%

\vskip .7 cm
\centerline{\autore
Giovanni Alberti, Guy Bouchitt\'e, Gianni Dal Maso}

%
%

\vskip .7 cm
\centerline{
\vbox{\hsize 10.3 cm\baselineskip 10 pt\parindent 0 pt
\eightrm 
{\eightbf Abstract.} 
In this paper we present a minimality criterion
for the Mumford-Shah functional, and more generally for non convex
variational integrals on {\eightti SBV} which couple a surface and a bulk
term. This method provides short and easy proofs for several minimality
results.
\smallskip
{\eightti Keywords}: Mumford-Shah functional, free discontinuity
problems, special functions of bounded variation, necessary condition
for minimality, calibrations, minimal partitions, gradient flow.
}}

%
%

\Section{\chp{s1}.~Introduction}
The Mumford-Shah functional was introduced in [\rf{MS1}] and [\rf{MS2}]
within the context of a variational approach to image segmentation
problems (cf.\ [\rf{MS2}] and [\rf{Mo-Sol}], Chapter 4).
In dimension $n$ it can be written as follows
$$
F(u)
:=\int_{\O\setminus Su} | \nabla u|^2 \,dx
        + \alpha\H^{n-1}(Su)
        + \beta\int_\O (u-g)^2 \,dx\ ,
\eqno (\frm{1.1})
$$
where $\O$ is a bounded regular domain in $\R^n$, $g: \O\to [0,1]$
is a given function (input grey level), $\alpha$ and $\beta$ are
positive (tuning) parameters, $\H^{n-1}$ is the $(n-1)$-dimensional
Hausdorff measure (that is, the usual $(n-1)$-dimensional area in case
of subsets of regular hypersurfaces, the length in the most relevant
case $n=2$). The unknown function $u: \O\to\R$ is regular (say, of
class $C^1$) out of a closed singular set $Su$, whose shape and
location are not prescribed.
Thus minimizing $F$ means optimizing the function $u$ and
the singular set $Su$. Indeed, in the original formulation
only the planar case $n=2$ is considered, and the singular set
is explicitly viewed as an independent variable, $u$ being
smooth on the complement of this set.

While existence results for minimizers of $F$ in dimension two
could be proved within the original framework (cf.\ [\rf{MS2}], 
[\rf{DMS}],
or [\rf{Mo-Sol}], Chapter 15), in arbitrary dimension they were
first obtained by a different approach outlined by E. De Giorgi
(cf.\ [\rf{DGA}]). More precisely, $F$ can be defined
for every function $u$ in the space $SBV(\O)$ of special functions
with bounded variation (see Section~\rf{s2} for more details, 
or [\rf{AFP}], Chapter 4; this space can be regarded
as a sort of the completion of piecewise regular functions),
upon which it is lower semicontinuous and coercive with respect
to the $L^1$ topology (cf.\ [\rf{Amb-ARMA}], or [\rf{AFP}, Chapter 6).
This immediately yields minimizers of $F$ within this space,
which can be proved to be minimizers of the original functional by
suitable regularity theorems (cf.\ [\rf{DGCL}], see also [\rf{AFP}],
Chapter 6, for further regularity results).
Furthermore, the lower semicontinuity and coercivity results
in [\rf{Amb-ARMA}] apply to a larger class of functionals
coupling bulk and surface contributions,
among which $F$ should be regarded as a prototypical example.
Consequently, the $SBV$ setting has been used to model a wide range
of problems, from image segmentation, to fractures in brittle
materials, to nematic liquid  crystals (see [\rf{AFP}], Section 4.6,
for a survey).

On a mathematical level, one of the most relevant features of the
Mumford-Shah functional is a deep lack of convexity. Hence, not only
minimizers may be not unique, but ``identifying'' them is by no means
an easy task, also in terms of efficient algorithms. Clearly,
every minimizer $u$ of $F$ must satisfy certain equilibrium conditions
-- an Euler-Lagrange equation of a sort -- which can be obtained
by considering different types of infinitesimal variations 
(see for instance [\rf{MS2}] or [\rf{AFP}], Section 7.4).
Among these we mention the following:
$u$ satisfies $\Delta u=\beta(u-g)$ in the complement
of the singular set $Su$, the normal derivative of $u$
on $Su$ must vanish (where $Su$ is a regular surface),
while the mean curvature of $Su$ multiplied by $\alpha$ is equal
to the difference of the energy densities $|\nabla u|^2+\beta(u-g)^2$
on the two sides of $Su$ (the first two correspond
respectively to the Euler-Lagrange equation for the convex
functional $\int_{\O\setminus K} [|\nabla u|^2+\beta(u-g)^2]\,dx$
and the associated Neumann boundary condition on $\bd K$;
together they are equivalent to
minimality with {\it prescribed\/} singular set $Su=K$).
However, due to the lack of convexity of $F$, these
conditions do not imply minimality -- not even local minimality.

In this paper we propose a sufficient condition for minimality,
and describe some applications.

We remark that the problem of finding sufficient conditions for 
minimality already makes sense for a simplified version of $F$ 
which occurs in the theory of interior regularity for minimizers 
of $F$, and is obtained by dropping the lower order term in 
(\rf{1.1}) and setting for simplicity $\alpha=1$, that is, 
$$
F_0(u)
:=\int_{\O\setminus Su} \hskip-7pt | \nabla u|^2\,dx  
+ \H^{n-1}(Su) \ .
\eqno (\frm{1.2})
$$
For the time being we focus on minimizers of $F_0$ with 
prescribed boundary values (in short, {\it Dirichlet minimizers\/}), 
and describe the basic idea of this paper, without dwelling on details. 
More precise definitions and accurate statements will be given in
Section~\rf{s3}.

Assuming that $u$ and $Su$ are sufficiently regular, let $u^+$ 
and $u^-$ denote the limits of $u$ on the two sides of $Su$, 
so that $u^+>u^-$, and let $\nu_u$ be the unit normal to $Su$ 
pointing from the side of $u^-$ to that of $u^+$; the 
{\it complete graph\/} of $u$ is the boundary of the subgraph 
of $u$ (the set of all points $(x,t)\in\O\times\R$ such that 
$t\le u(x)$), oriented by the inner normal $\nus{\GM u}$. 
Thus $\GM u$ consists of the union of the usual graph of $u$ 
and an additional part given by all segments with endpoints 
$(x,u^-(x))$ and $(x,u^+(x))$, with $x$ ranging in
$Su$; see Figure 1 (cf.\ also Remark \rf{s2.3}).

\figureps{1}{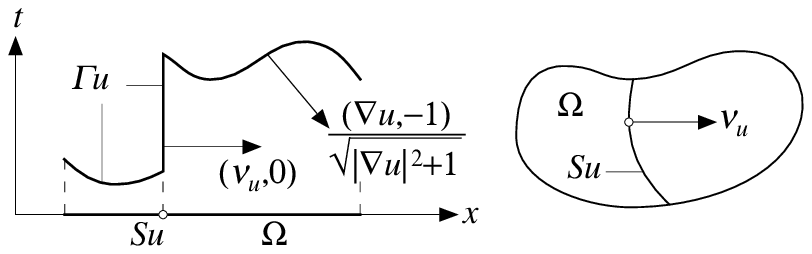}

We consider now the vectorfields $\phi=(\phi^x,\phi^t)$ on
$\O\times\R$ such that, for every function $u$, $F_0(u)$ is
larger than or equal to the flux of $\phi$ through $\GM u$, that is
$$
F_0(u)\ge\int_{\GM u} \phi\cdot \nus{\GM u} \, d\H^n\ .
\eqno (\frm{1.3})
$$
Taking into account (\rf{1.2}) and the fact that
the flux of $\phi$ through $\GM u$ is given by (cf.\ Figure 1)
$$
\int_\O \big[ \phi^x(x,u) \cdot \nabla u - \phi^t(x,u) \big]
    \, dx
+\int_{Su} \Big[
    \int_{u^-}^{u^+} \phi^x(x,t) \, dt
  \Big]  \cdot \nu_u \, d\H^{n-1} \ ,
\eqno (\frm{1.4})
$$
we see that inequality (\rf{1.3}) is satisfied if
$$
\phi^x(x,u) \cdot \nabla u - \phi^t(x,u) \le  |\nabla u|^2
\quad\hbox{for every $x\in\O\setminus Su$}
$$
and
$$
\Big| \int_{u^-}^{u^+} \phi^x(x,t) \, dt \Big| \le 1
\quad\hbox{for every $x\in Su$.}
$$
The first inequality holds for every $u$
if $\phi$ satisfies
$\phi^x (x,t) \cdot \xi - \phi^t(x,t)\le |\xi|^2$
for every point $(x,t)$ and for every vector $\xi$,
which can be equivalently restated as

\smallskip
\itemm{(a)}
{$|\phi^x(x,t)|^2 \le 4 \phi^t(x,t)$ for $x\in\O$, $t\in\R$,}

\smallskip\noindent
while the second inequality is satisfied if

\smallskip
\itemm{(b)}
{$\dps\Big| \int_{t_1}^{t_2} \phi^x(x,t) \, dt \Big|\le 1$
for $x\in\O$, $t_1,t_2\in\R$.}

\smallskip\noindent
Moreover, one easily checks that equality holds in (\rf{1.3})
for a particular $u$ if (and only if)

\smallskip
\itemm{(a')}
{$\phi^x(x,u(x))=2\nabla u(x)$ and
$\phi^t(x,u(x))=|\nabla u(x)|^2$ for
$x\in\O\setminus Su$,}

\itemm{(b')}
{$\dps\int_{u^-(x)}^{u^+(x)} \phi^x(x,t) \, dt=\nu_u(x)$ 
for $x\in Su$.}

\medskip
Let now be given a function $u$, and assume that there exists 
a vectorfield $\phi$ which is {\it divergence-free\/} and satisfies
assumptions (a), (b), (a'),
and (b') above. Then for every function $v$ which agrees with $u$
on the boundary of $\O$ we have
$$
F_0(v)
\ge \int_{\GM v} \phi \cdot \nus{\GM v} \, d\H^n
=   \int_{\GM u} \phi \cdot \nus{\GM u} \, d\H^n
=F_0(u) \ ,
\eqno(\frm{1.5})
$$
where the first equality follows from the divergence theorem,
since $\phi$ is divergence-free and $\GM u$ and $\GM v$ have the same
boundary. Hence the existence of such a vectorfield $\phi$
implies the minimality of $u$.
In Section~\rf{s3} we give a more precise version of this
result, and extend it to minimizers of $F$ and other
functionals (cf.\ Theorems \rf{s3.2.0} and \rf{s03.15}).
In particular, for minimal partitions we recover 
the principle of {\it paired calibrations\/} introduced in 
[\rf{ML1}] and [\rf{Bra1}] (in fact, a slight generalization 
of it -- cf.\ Theorem \rf{s03.16} and Remark \rf{s03.17}).

We call $\phi$ a {\it calibration\/} (for $u$) by analogy with the
corresponding theory for minimal surfaces: in that setting,
a vectorfield $\phi$ is said to calibrate an oriented hypersurface
$S$ (with boundary) if it agrees on $S$ with the normal
vectorfield, is divergence-free, and satisfies $|\phi|\le 1$
everywhere; the existence of a calibration implies that $S$
minimizes the area among all oriented hypersurfaces
with the same boundary, and the proof is just one line, as (\rf{1.5})
above (we refer the reader to [\rf{Mor}] for detailed references
and a review of many results). Calibrations have also been defined
for general integrals on oriented surfaces (or currents)
of any dimension and codimension (cf.\ [\rf{Fed-calib}],
notice that in this general framework they are closed
differential form rather than divergence-free vectorfields).
Moreover, since every variational integral of the form 
$\int_\O f(x,u,\nabla u)\,dx$ with $u:\O\to\R^k$ and $f$ polyconvex 
in $\nabla u$ can be re-written as a one-homogeneous convex integral 
on the graph of $u$, the theory of calibrations can be adapted also 
to these problems, yielding the equivalent of null-lagrangians.

However, this geometric interpretation does not cover 
free-discontinuity problems, because functionals of
Mumford-Shah type cannot be written as integral over 
the (complete) graph of $u$. Indeed, the novelty of our 
approach consists in introducing suitable non-local constraints
(namely, condition (b)) to define the class of admissible calibrations. 

\smallskip
Once the general principle of calibrations is given, two fundamental
questions arise: does every minimizer admit a calibration?
And how can we recover it?
In other words, we have given a sufficient condition for minimality,
but we do not know if it is actually fulfilled by any minimizer,
and, which is even more relevant in applications, we do not know
how to verify it, that is, how to construct a calibration.

To discuss the first question, we begin by recalling the situation
of minimal surfaces. The basic idea behind calibrations is that the
area functional admits a natural extension from regular oriented
surfaces (with fixed boundary) to all normal currents
(with the same boundary), which is nothing else but the mass.
But now we have a convex functional on an affine space,
and hence minimum points are exactly characterized by the fact
that the subdifferential of the functional contains the zero
element; in the specific case, this turns out to be equivalent
to the existence of a calibration with Borel coefficients
(cf.\ [\rf{Fed-calib}], Proposition 4.10(3)).
Thus a given surface admits a calibration of a sort if and only
if it minimizes the area (the mass) among all {\it normal\/} currents
with same boundary. However, in codimension (and dimension)
larger than one, minimizing the area in the class of surfaces,
or even integral currents, does not necessarily imply minimizing
the area in the class of normal currents. In other words,
the infima of the area functional on the two classes
may not coincide: counterexamples were given, for instance,
in [\rf{Mor2}] and [\rf{White}].

The situation of the Mumford-Shah functional is somehow similar.
Let $\G_0$ be the class of all vectorfields $\phi$ on $\O\times\R$
which satisfy assumptions (a) and (b) above;
since for every function $u$ we can find $\phi\in\G_0$
which also satisfies assumptions (a') and (b'),
we have the identity
$$
F_0(u)
=\sup_{\phi\in\G_0} \int_{\GM u} \phi \cdot \nus{\GM u} \, d\H^n
=\sup_{\phi\in\G_0} \int_{\O\times\R} \phi \cdot D1_u\ ,
\eqno (\frm{1.6})
$$
where $D1_u$ is the (distributional) derivative of the
the characteristic function $1_u$ of the subgraph of $u$.
Thus $F_0(u)=G(1_u)$, where $G(v)$ is defined as the supremum 
of $\int_{\O\times\R} \phi \cdot Dv$ over all $\phi\in\G_0$ for 
every function $v\in BV\loc(\O\times\R)$. 
Since $G$ is a convex functional on the affine space
$X_u$ of all functions $v$ which agree with $1_u$ on the boundary,
$1_u$ minimizes $G$ if and only if the subdifferential
of $G$ at $1_u$ contains the zero element, condition which is 
roughly equivalent to the existence of a calibration for $u$.
It was proved in [\rf{Cha}] for the one-dimensional case $n=1$ 
that the function $1_u$ minimizes $G$ whenever $u$ minimizes $F_0$;
in other words, the infimum of $G$ on $X_u$ agrees with the infimum 
of $F_0$ (among all functions with same boundary values as $u$),
thus proving that every minimizer of $F_0$ can be calibrated
(in some sense).
It is not known if the same result holds in higher dimension.

The discussion of this point may benefit from being set in a more 
abstract framework. We can indeed summarize the calibration method
as follows: given a nonconvex function(al) $F$ on a certain space
$X$, we embed $X$ into a convex space $\tilde X$ so that 
$F$ agrees with a convex function $\tilde F$ on $\tilde X$. 
Thus the minimality on a point $x\in X$ with respect to $\tilde F$ 
can be characterized via the subdifferential of $\tilde F$ at $x$; 
this yields a sufficient condition for the minimality with respect 
to $F$, conditions that is  also necessary when (and only when) the 
infima of $\tilde F$ on $X$  and $\tilde X$ agree. 
As pointed out in [\rf{BC}], 
it is always possible to construct the convex space $\tilde X$ 
and the function $\tilde F$ so that the two infima agree, 
but the problem is that checking the condition of 
minimality derived from this abstract construction may be neither
simpler nor more feasible than a direct verification of minimality.

\smallskip
Passing to the second question mentioned above, we notice that an
abstract result which guarantees the existence of calibrations
would nevertheless provide no solution to the problem of construction. 
As a matter of fact, we do not know of any general method
of construction, not even for minimal surfaces.
Instead, we have collected in Sections~\rf{s4} and~\rf{s5} many examples
of calibrations for $F_0$ and $F$, and gathered some helpful
remarks and observations.
Despite a lack of a general recipe, the calibration method provides 
short and easy proofs of some natural minimality result; among these 
we recall the following ones:

\itemm{(1)}
{every harmonic function is a Dirichlet minimizer of
$\int_\O |\nabla u|^2 \, dx$, and is also a Dirichlet minimizer 
of $F_0$ when the gradient is sufficiently small (Paragraph \rf{e9});}

\itemm{(2)}
{a function which is constant on each element of a minimal partition
of the domain is a Dirichlet minimizer of $F_0$ when the
values are sufficiently far apart from each other;
in particular this applies to the so-called triple junction
(Paragraphs \rf{e11} and \rf{e13});}

\itemm{(3)}
{every solution of the equation $-\Delta u +\beta(u-g)=0$
with Neumann boundary conditions is a minimizer of
$\int _\O [|\nabla u|^2 +\beta(u-g)^2]\,dx$, and also of
$F$ for large $\beta$  (Paragraph \rf{e18});}

\itemm{(4)}
{if $g$ is the characteristic function of a regular set, then 
$u:=g$ minimizes $F$ for large $\beta$ (Paragraph \rf{e19}).}

\smallskip\noindent
Notice that (3) and (4) give a strong indication that for
initial data with smooth singular sets
the gradient flow associated with $F_0$ in the $L^2$ metric
(which can be defined via time discretization, cf. [\rf{Gob2}]) 
leaves the singular set still, at least for small times,
and agrees with the heat flow elsewhere (cf. [\rf{Gob1}]
and Remarks \rf{e25}, \rf{e26}, \rf{e27}).

\smallskip
We conclude this introduction with a few remarks.
A long-standing conjecture on the Mumford-Shah functional is that 
the ``cracktip'', namely the function $u$ on the plane given in polar 
coordinates by $u:=\sqrt{2\rho/\pi}\,\sin(\theta/2)$, with 
$-\pi<\theta\le\pi$, (or, equivalently, the imaginary part of 
$\sqrt{2z/\pi}$ with a cut along the negative real axis) 
is a Dirichlet minimizer of the homogeneous
Mumford-Shah functional $F_0$ on every bounded open subset of $\R^2$.
This conjecture has been recently proved in [\rf{Bon-Dav}],
but so far no calibration has been found for this minimizer.

As shown in [\rf{Bra3}] in the case of minimal partitions, 
the calibration method can be numerically implemented to get 
rigorous lower bounds for the value of minima, to disprove 
the minimality of a given configuration, or viceversa to get 
an idea of what the calibration (if it exists) should look like. 
It is not clear if an efficient numerical implementation for
the Mumford-Shah functional is feasible. 

The theory presented in this paper is limited 
to scalar functions. To develop a similar theory for functionals 
on $\R^k$-valued maps $u$, one should replace 
divergence-free vectorfields by closed $n$-forms on $\O\times\R^k$,
which act on the graphs of the maps $u$, viewed as
(suitably defined) $n$-surfaces in $\O\times\R^k$.

\medskip
This paper is organized as follows. In Section \rf{s2}
we recall the basic notation about finite perimeter sets
and the space $SBV$, which is indeed the natural setting
for our theory (however, under most regards the unfamiliar 
reader can just replace the word $SBV$ with ``smooth out 
of a piecewise smooth singular set'', and ``finite perimeter'' 
with ``piecewise smooth boundary'').
In Section \rf{s3} we expand the idea outlined above and develop 
the theory of calibrations for minimizers of $F_0$ and $F$,
with or without prescribed boundary values. 
Then we discuss the extension to more general functionals, 
and the connection with paired calibrations.
Sections \rf{s4} and \rf{s5} are devoted to 
examples and applications.
Finally, the appendix contains the proofs of some technical
results stated in Section \rf{s2}.

Some of the results contained in this paper were announced in
[\rf{ABD}] and [\rf{3ecm}]. Further applications can be found in
[\rf{DM-M-M}] and
[\rf{M-M}].

\Parag{Acknowledgements}
The first and third authors have been partially supported 
by MURST through the projects ``Equazioni Differenziali e Calcolo 
delle Variazioni (1997)'' and ``Calcolo delle Variazioni (2000)''.
This research was initiated while the first author was visiting 
the University of Toulon, and subsequently developed during a stay 
at the Max Planck Institute for Mathematics in the Sciences in Leipzig.
Several people contributed, with discussions and remarks, to the final 
shape of this paper; among them, we would like to thank in particular 
Antonin Chambolle and Massimo Gobbino.


\Section{\chp{s2}.~Notation and preliminary results}
Throughout this paper, sets and functions are always assumed
to be Borel measurable, and we {\it do not} identify functions
which agree almost everywhere.
A {\it vectorfield\/} on a subset $E$ of $\R^n$ is any map from $E$ 
into $\R^n$. The divergence of a vectorfield is always intended in 
the sense of distributions (in the interior of its domain); in particular
we say that $\phi$ is divergence-free to mean that its divergence
vanishes.  

The letter $\O$ denotes a (possibly unbounded) open subset of $\R^n$,
$\SS^{n-1}$ is the unit sphere in $\R^n$,
$\H^k$ stands for the $k$-dimensional Hausdorff measure
(which agrees with the usual $k$-dimensional volume
on every regular surface of dimension $k$), and $\L^n$ is
the $n$-dimensional Lebesgue measure;
when the integration is done with respect to
Lebesgue measure, we always write $dx$ instead of $d\L^n$.
The restriction of any Borel measure
$\mu$ to a set $E$ is denoted by $\mu\LL E$, while
$g\cdot\mu$ is the (vector) measure canonically associated with any
$\mu$-summable (vector) function $g$. Thus $\mu\LL
E=1_E\cdot \mu$ where $1_E$ is the characteristic function of $E$,
namely the function which is equal to $1$ on $E$ and to $0$ elsewhere.
$\|\cdot\|_p$ denotes the norm in the space $L^p$.

A (vector) function $f$ on $\R^n$ has {\it approximate limit\/}
$a$ at $x$, and we write $\aplim\limits_{y\to x} f(y)=a$, if 
$$
\lim_{r\to 0} {1\over r^n} \int_{B(x,r)} |f(y)-a| \, dy =0 \ .
\eqno(\frm{a1})
$$
The same definition applies to functions defined on a subset $E$ 
of $\R^n$, provided that $B(x,r)$ is replaced by $E\cap B(x,r)$. 
Notice that this definition slightly differs from the more usual 
one, which is expressed in term of the density at $x$ of the 
pre-images of neighbourhoods of $a$ (cf.\ [\rf{EG}], Section 1.7.2, 
and [\rf{AFP}], Definition 3.63); however, these two notions agree 
for locally bounded functions. 

We recall now some notation and basic facts about finite perimeter sets,
$BV$ and $SBV$ functions; for a more precise definitions and a detailed
account of the results we refer to [\rf{AFP}], Chapters 3 and 4 
(for the theory of $BV$ functions, see also [\rf{EG}], Chapter 5).

A real function $u$ on an open set $\O\subset\R^n$ has {\it bounded
variation\/} in $\O$, and we write $u\in BV(\O)$, if it belongs to
$L^1(\O)$ and the distributional
derivative $Du$ is (represented by) a bounded vector measure on $\O$.
The integral of an $\R^n$-valued function $f$ with respect to the 
measure derivative $Du$ is denoted by $\int_\O f\cdot Du$.
Every $u\in BV(\O)$ is almost everywhere differentiable
in the approximate (or measure theoretic) sense, and the corresponding
{\it approximate gradient\/} $\nabla u$ agrees with the density of the 
measure $Du$ with respect to Lebesgue measure.
If $\O$ has Lipschitz boundary, then $u$ admits a {\it trace\/} on 
the boundary (in the approximate sense), which we still denote by $u$, 
and which belongs to $L^1\loc(\bd\O,\H^{n-1})$.

The {\it singular set\/} $Su$ is the set of all points where $u$ admits 
no approximate limit; $Su$ has Hausdorff dimension (less than or) equal to
$n-1$, and more precisely it is {\it rectifiable\/} (of dimension $n-1$),
which means that it can be covered, up to an
$\H^{n-1}$-negligible subset, by countably many hypersurfaces of class
$C^1$ (these sets are sometimes called 
``countably $(\H^{n-1}, n-1)$-rectifiable'', cf.\ [\rf{Fed-GMT}], Section
3.2.14). We recall that, for every rectifiable set, the {\it approximate
unit  normal\/}, and the corresponding approximate tangent hyperplane, 
are well-defined at $\H^{n-1}$-almost every point (and do not depend
on the choice of the covering). The approximate unit normal 
to $Su$ at $x$ is denoted by $\nu_u(x)$.
At $\H^{n-1}$-almost every $x\in Su$ there exist the approximate 
limits $u^+(x)$ and $u^-(x)$ of $u$ on the two sides of $Su$ 
(more precisely, the approximate limits of the restriction of $u$ 
to the two half-spaces defined by the approximate tangent 
hyperplane at $x$), and we arrange so that $u^+(x) > u^-(x)$, 
and $\nu_u(x)$ is pointing  from the side of $u^-(x)$ to the 
one of $u^+(x)$. 

The measure $Du$ can be canonically decomposed as the sum of 
three mutually orthogonal measures: the {\it Lebesgue part\/} 
$\nabla u\cdot \L^n$, the {\it jump part\/} 
$(u^+ - u^-) \, \nu_u \cdot\H^{n-1} \LL Su$, and a remainder, 
called {\it Cantor part\/}, which is singular but does not charge any
$\H^{n-1}$-finite set.

The space $SBV(\O)$ of {\it special functions of bounded variation\/}
is given by all functions
$u\in BV(\O)$ for which the Cantor part of the derivative vanishes,
i.e.,
$$
Du=\nabla u \cdot \L^n + (u^+ - u^-) \, \nu_u \cdot \H^{n-1} \LL Su \ .
\eqno (\frm{2.1})
$$

A subset $E$ of $\O$ has {\it finite perimeter\/} (in $\O$) if the 
distributional derivative $D1_E$ of its characteristic function  
is a bounded vector measure on $\O$. 
In this case, the {\it measure theoretic boundary\/}
$\bd_* E$ is the singular set of $1_E$,
while the inner normal $\nus{\bd_*E}$ is
the associated normal vectorfield $\nu_{1_E}$.
Both the Lebesgue and the Cantor parts of $D1_E$ vanish,
i.e.,
$$
D1_E=\nus{\bd_* E} \cdot\, \H^{n-1}\LL\bd_*E \ .
\eqno (\frm{2.2})
$$

\subsection{A general form of the divergence theorem}
When looking for calibrations for a given function, 
it is often convenient to consider also vectorfields which
are not regular. In doing so, however, we face some 
technical difficulties. 

First of all, the first identity in (\rf{1.5}) depends 
on the divergence theorem, and may not hold when $\phi$ is 
divergence-free (in the sense of distribution) but not continuous, 
because the flux of such a vectorfield through a given surface
is not well-defined. To solve this problem, we must assume 
a certain regularity in $\phi$.

\Def{\thm{s2.0.0}}
{\it
We say that a vectorfield $\phi$ on a subset $E$ of $\R^n$
is {\rm approximately regular} if it is  bounded, 
and for every Lipschitz hypersurface $M$ in $\R^n$
there holds
$$
\aplim_{y\to x}\big[ \phi(y) \cdot \nus{M}(x) \big]
= \phi(x) \cdot \nus{M}(x)
\quad\hbox{for $\H^{n-1}$-a.e.\ $x\in M\cap E$,}
\eqno (\frm{2.2.1})
$$
where $\nus{M}(x)$ denotes the (unit) normal to $M$ at $x$.
}

\Rem{\thm{s2.0.1}}
If $\phi$ is approximately regular, then (\rf{2.2.1}) can be 
extended to every rectifiable set $M$, 
$\nus{M}$ being now understood in the approximate sense.
If $\phi$ admits traces $\phi^+$ and $\phi^-$ on the two sides of 
$M$ (defined in the same way of the traces $u^+$ and $u^-$), then 
(\rf{2.2.1}) is equivalent to the  compatibility condition
$$
\phi \cdot \nus{M} = \phi^+ \cdot \nus{M} = \phi^- \cdot \nus{M}
\quad\hbox{$\H^{n-1}$-a.e.\ in $M\cap E$,}
\eqno (\frm{2.2.2})
$$
which links the pointwise values of $\phi$ on $M$ with the values of
the traces.

In particular, if $\phi$ is (approximately) continuous
$\H^{n-1}$-almost everywhere on $E$, then it is approximately
regular. Similarly, if $\phi$ is (approximately) continuous on the 
complement of a rectifiable set $S$, then $\phi$ is 
approximately regular if and only if (\rf{2.2.1}) holds for $M:=S$.

\Rem{\thm{r1}}
If $\phi$ has the special form $\phi:=(0,\ldots, 0, \psi)$
where $\psi=\psi(x_1,\dots,x_n)$ is a bounded real function 
which is continuous in the variable $x_n$, 
then $\phi$ is approximately regular. 
Take indeed a Lipschitz surface $M$, and let $M_0$ be the subset of
all points $x\in M$ such that the $n$-th component of $\nus{M}(x)$ 
vanishes. 
Then equality (\rf{2.2.1}) obviously holds for all $x\in M_0$. 
To prove that it also holds for $\H^{n-1}$-a.e.\ point in 
$M\setminus M_0$ it suffices to notice that $\psi$ is approximately
continuous in the complement of a set of type $N\times\R$, where $N$ 
is a $\L^{n-1}$-negligible  subset of $\R^{n-1}$. 
Thus $\phi$ is approximately continuous at all points  except those 
in $N\times\R$, which form an $\H^{n-1}$-negligible subset of
$M\setminus M_0$ by the area formula (see [\rf{Fed-GMT}], Theorem 
3.2.22, or, for a slightly less general statement, [\rf{AFP}], 
Theorem 2.91). 

\medskip
We can now state a refined version of the classical divergence theorem 
(the proof is postponed to the appendix). 

\Lemma{\thm{s2.0.2}}
{\it
Let $\O$ be an open set in $\R^n$ with Lipschitz boundary, 
$\phi$ an approximately regular vectorfield on $\ove\O$, and 
$u$ a function in $BV(\O)$. Assume moreover that 
$\div\phi\in L^\infty(\O)$ and $u\phi\in L^1(\bd\O,\H^{n-1})$. Then 
$$
\int_\O \phi \cdot Du
= -\int_\O  u \,\div\phi \, dx
-\int_{\bd\O}  u \,\phi\cdot\nus{\bd\O} \, d\H^{n-1} \ ,
\eqno (\frm{2.2.3})
$$
where $\nus{\bd\O}$ is the inner unit normal to $\bd\O$, and in the last 
integral $u$ stands for the trace of $u$ on $\bd\O$.
}

\medskip 
Notice that the condition  $u\phi\in L^1(\bd\O,\H^{n-1})$ is always 
satisfied when $\bd\O$ is bounded, because in this case the trace of $u$
belongs to $L^1(\bd\O,\H^{n-1})$.

\Rem{\thm{s2.0.3}}
As shown in [\rf{Anz}], for a Borel vectorfield $\phi$ with divergence 
in $L^1$ it is possible to give a functional definition of the trace 
of the normal component of $\phi$ on any Lipschitz surface $M$, and 
precisely as a continuous extension of the trace operator for regular
vectorfields. 
When $M$ is the boundary (of some set $\O$), this notion of
trace automatically  satisfies the divergence theorem (in the sense that
formula (\rf{2.2.3}) holds for every function $u$ of class $C^1_c$), but
since modifying $\phi$ in a Lebesgue-negligible set affects neither this
trace nor the distributional divergence, one can easily produce examples
where the trace does not agree with  the normal component of $\phi$ on $M$.
Thus Lemma \rf{s2.0.2} shows implicitly  that this is never the case when
$\phi$ is approximately regular. One may wonder if every vectorfield with
divergence in $L^\infty$ agrees, up to a Lebesgue-negligible set, with an
approximately regular one.  Unfortunately the answer is negative, even for
divergence-free vectorfields (cf.\ [\rf{Anz}], example after Proposition
2.1).

\medskip
Going back to the first identity in (\rf{1.5}), we remark that 
verifying that a vectorfield $\phi$ is divergence-free is relatively
easy when $\phi$ is of class $C^1$ because the distributional
divergence agrees with the classical one, which can be explicitly 
computed. If $\phi$ is piecewise $C^1$, the task is slightly
more difficult, and can be carried out in many concrete cases 
(see Sections \rf{s4} and \rf{s5}) with the help of the following lemma
(the proof is postponed to the appendix).

\Lemma{\thm{s2.7}}
{\it
Let $\phi$ be a bounded vectorfield on an open set $\O\subset\R^n$, 
and assume that there exist a closed set $S$
and a function $f\in L^1\loc(\O)$ such that $\div\phi=f$ in the sense 
of distributions on $\O\setminus S$.
Then the identity $\div\phi=f$ holds also on $\O$ if
$S$ can be written as $S:=S_0\cup S_1$, with 
$S_0$ an $\H^{n-1}$-negligible closed set and $S_1$ a 
(possibly disconnected) Lipschitz hypersurface, and
$\phi$ satisfies {\rm (\rf{2.2.1})} for $M:=S_1$ and $E:=\O$. 
}

\Rem{\thm{s2.8}}
The point of this lemma is roughly the following: 
since the divergence is a first order differential operator, $\div\phi$ 
cannot ``charge'' any set of codimension larger than 1, and therefore 
$S_0$ can be safely removed. On the other hand, the part of $\div\phi$
supported  on the hypersurface $M$ is given by the difference of the 
traces (whenever defined) of the normal components of $\phi$ on the two
sides of $M$, which happens to vanish if (\rf{2.2.1}) holds, and then 
we are allowed to neglect $S_1$ too.

\Rem{\thm{s2.9}}
Lemma \rf{s2.7} will be often applied in the following forms.

(a)~Suppose that $\phi$ is a bounded vectorfield on $\ove\O$, continuous 
on $\ove \O\setminus (S_0\cup S_1)$, and divergence-free on
$\O\setminus (S_0\cup S_1)$, with $S_0$, and $S_1$ given as above. 
If $\phi$ satisfies (\rf{2.2.1}) with $M=S_1$, then $\phi$ is
approximately regular on $\ove\O$ and divergence-free on $\O$
(cf. Remark \rf{s2.0.1}).

(b)~Let be given, for $j=1,\dots,m$, pairwise disjoint 
Lipschitz open sets $\O_j$ whose closures cover $\ove\O$, 
and approximately regular, divergence-free vectorfields $\phi_j$ 
on $\ove\O_j$. Let $\phi$ be any vectorfield on $\ove\O$ which agrees
at any point with one of the $\phi_i$ (hence it 
is uniquely determined at least on the union of all $\O_i$).
Then $\phi$ is approximately regular and divergence-free
provided that the vectorfields $\phi_i$ satisfy the compatibility 
conditions
$$
\phi_i\cdot\nus{\bd\O_i} =\phi_j\cdot\nus{\bd\O_j}
\quad\hbox{$\H^{n-1}$-a.e.\ on $\bd\O_i\cap\bd\O_j$,}
$$
which are equivalent to the compatibility condition (\rf{2.2.2}) 
for $\phi$.

\subsection{The complete graph of an $SBV$ function}
We fix now some notation and state some results which are more 
specific to this paper.
In the following $\O$ is a fixed bounded open subset of $\R^n$ 
with Lipschitz boundary, and $\nus{\bd\O}$ is its inner unit normal.
The letter $x$ usually denotes the variable in $\O$ (or $\R^n$),
while $t$ is the variable in $\R$; $U$ is an open subset of $\O\times\R$,
with Lipschitz boundary whose closure can be written as
$$
\ove U:=\big\{ (x,t)\in\ove \O\times\R : \ \tau_1(x)\le t\le 
\tau_2(x) \big\} \ ,
\eqno (\frm{2.2.4})
$$
where the functions $\tau_1, \tau_2:\overline\O\to[-\infty,+\infty]$
satisfy $\tau_1<\tau_2$.
The letter $\phi$ denotes a bounded vectorfield defined on (a subset of) 
$\R^n\times\R$, with components $\phi^x\in\R^n$ and $\phi^t\in\R$. 
Notice that $\div\phi=\div_x\phi^x+\bd_t\phi^t$, where $\div_x$ is 
the (distributional) divergence with respect to $x$ 
and $\bd_t$ the (distributional) derivative  with respect to~$t$.

\Def{\thm{s2.1}}
{\it
For every function $u\in BV(\O)$, let $1_u$ be the characteristic
function of the {\rm subgraph} of $u$ in $\O\times\R$, namely 
$1_u(x,t):=1$ for $t\le u(x)$ and $1_u(x,t):=0$ for $t> u(x)$.
The {\rm complete graph} of $u$, denoted by $\GM u$, is the measure
theoretic boundary of the subgraph of $u$, i.e., the singular set 
of $1_u$.
}

\medskip
Since the subgraph of $u$ has finite perimeter in $\O\times\R$ (see,
e.g., [\rf{Mir}], Proposition 1.4), the definition of $\GM u$ is 
well-posed.
Moreover $D1_u=\nus{\GM u} \cdot \H^n \LL \GM u$ (cf.\ (\rf{2.2})),
where $\nus{\GM u}$ is the inner unit normal of the subgraph of $u$. 
Therefore the flux through $\GM u$ of any vectorfield
$\phi$ on $\O\times\R$ is given by the integration of $\phi$ with
respect to the vector measure $D1_u$, that is
$$
\int_{\GM u} \phi \cdot \nus{\GM u} \, d\H^n
=\int_{\O\times\R} \phi \cdot D1_u \ .
\eqno (\frm{2.3})
$$
An alternative way to compute this flux is given by the following lemma
(the proof is postponed to the appendix).

\Lemma{\thm{s2.2}}
{\it
Let $u$ be a function in $SBV(\O)$ and let $\phi$ be a
vectorfield defined at least on $\GM u$. Then
$$
\eqalign{
\int_{\GM u} \phi \cdot \nus{\GM u} \, d\H^n 
=&  \int_\O \big[ \phi^x(x,u) \cdot \nabla u - \phi^t(x,u) \big] \, dx \cr
 &+ \int_{Su} \Big[ \int_{u^-}^{u^+} \phi^x(x,t) \, dt \Big] \cdot \nu_u
    \, d\H^{n-1}\ , \cr
}
\eqno (\frm{2.4})
$$
where $u$, $u^\pm$, $\nabla u$, and $\nu_u$ are always computed at $x$.
}

\Rem{\thm{s2.3}}
Formula (\rf{2.4}) corresponds to a decomposition of the derivative
of $1_u$, or, better, to a decomposition of the complete graph $\GM u$
as union (up to $\H^n$-negligible sets) of a ``regular'' part --
namely the set of all points $\big(x,u(x)\big)$ such that $u$ is
approximately continuous at $x$, and has approximate gradient $\nabla
u(x)$ -- and a ``vertical'' part -- namely the set of all points
$(x,t)$ with $x\in Su$ and $t\in (u^-(x), u^+(x))$. Note that
for a general $BV$ function there would be an additional subset of
$\GM u$, corresponding to the Cantor part of $Du$.

\medskip
The following version of the divergence theorem (cf.\ Lemma \rf{s2.0.2})
yields the first equality in (\rf{1.5}) (the proof is postponed to the 
appendix).

\Lemma{\thm{s2.6}}
{\it
Let be given two functions $u$ and $v$ in $BV(\O)$ whose complete graphs
lie in $\ove U$, and an approximately regular vectorfield $\phi$ on 
$\ove U$ which is divergence-free in $U$. Then
$$
\eqalign{
    \int_{\GM u} \phi \cdot \nus{\GM u} \, d\H^n
& - \int_{\GM v} \phi \cdot \nus{\GM v} \, d\H^n =\cr
& = \int_{\bd\O} \Big[ \int_u^v \phi^x(x,t) dt\Big]
    \cdot \nus{\bd\O} \, d\H^{n-1} \cr
}
\eqno (\frm{2.7})
$$
(where, in the last integral, $u$ and $v$ stand for the traces on $\bd\O$).
}

\Section{\chp{s3}.~Calibrations for free discontinuity problems}
In this section we introduce a calibration principle for a wide class of
free discontinuity problems, expanding the basic idea explained in
the introduction. We begin with the case of the
Mumford-Shah functional, with or without the lower order term,
and then we consider more general functionals, possibly with
discontinuous integrands, which include some interesting functionals
considered in minimal partition problems.

Throughout this section $\O$ is a bounded open subset of $\R^n$ with
Lipschitz boundary, $U$ is an open set with Lipschitz boundary 
contained in $\O\times\R$ which satisfies (\rf{2.2.4}); 
$u$ always denotes a function in $SBV(\O)$. 
The functionals $F(u)$ and $F_0(u)$ are given in (\rf{1.1}) and (\rf{1.2}),
respectively, where $Su$ and $\nabla u$ are now defined as in 
Section~\rf{s2}, $\alpha>0$ and $\beta\ge0$ are fixed constants, and $g$
belongs to $L^\infty(\O)$. Note that the functional $F_0$ is the particular
case of $F$ corresponding to $\alpha=1$ and $\beta=0$.
In the following definition we fix some terminology about minimizers 
of $F$ (which also applies to any other functional on $SBV$). 

\Def{\thm{s3.0}}
{\it
We say that a function $u$ is an {\rm (absolute) minimizer} of 
$F$ if $F(u)\le F(v)$ for all $v\in SBV(\O)$, while $u$ is a 
{\rm Dirichlet minimizer} if $F(u)\le F(v)$ for all $v\in SBV(\O)$ 
with same trace on $\bd\O$ as $u$ (that is, with same boundary 
values as $u$). We say that $u$ is a $\ove U$-minimizer 
if the complete graph of $u$ is contained in $\ove U$ and $F(u)\le F(v)$ 
for all $v\in SBV(\O)$ with complete graph contained in $\ove U$, 
while $u$ is a $\ove U$-Dirichlet minimizer if we add the requirement 
that the competing functions $v$ have the same boundary values as $u$.
}

\subsection{Calibrations for the Mumford-Shah functional}
We generalize the idea described in the introduction and provide
sufficient conditions for $\ove U$-minimality and $\ove U$-Dirichlet 
minimality with respect to the Mumford-Shah functional $F$ (or $F_0$).
We begin with the following key lemma.

\Lemma{\thm{s3.1.0}}
{\it
Let $F$ be defined as in {\rm(\rf{1.1})} for some $\alpha>0$, 
$\beta\ge 0$, and let $\phi$ be a vectorfield on $\ove U$ which 
satisfies the following  assumptions:

\smallskip
\itemm{(a)}
{$\phi^t(x,t)\ge {1\over 4} |\phi^x(x,t)|^2-\beta(t-g)^2$
for $\L^n$-a.e.\ $x\in\O$ and every
$t\in [\tau_1,\tau_2]$,}

\smallskip
\itemm{(b)}
{$\dps\Big| \int_{t_1}^{t_2} \phi^x(x,t) \, dt \Big|\le\alpha$
for $\H^{n-1}$-a.e.\ $x\in\O$ and every
$t_1,t_2\in [\tau_1,\tau_2]$,}

\smallskip\noindent
where the functions $\tau_1$ and $\tau_2$ are defined 
by (\rf{2.2.4}) and, like $g$, are computed at $x$.
Then for every $u$ such that $\GM u\subset\ove U$ we have
$$
F(u)\ge \int_{\GM u} \phi \cdot \nus{\GM u} \, d\H^n \ .
\eqno (\frm{3.1.0})
$$
Moreover, equality holds in {\rm (\rf{3.1.0})} for a given
$u$ if and only if

\smallskip
\itemm{(a')}
{$\phi^x(x,u)=2\nabla u$ and
$\phi^t(x,u)=|\nabla u|^2 - \beta(u-g)^2$ 
for $\L^n$-a.e.\ $x\in\O$,}

\itemm{(b')}
{$\dps\int_{u^-}^{u^+} \phi^x(x,t) \, dt=\alpha\,\nu_u$
for $\H^{n-1}$-a.e.\ $x\in Su$,}

\smallskip\noindent
where $u$, $u^\pm$, $\nabla u$, $\nu_u$, and $g$
are always computed at $x$.
}

\Pr
Take $u$ such that $\GM u\subset\ove U$.
We recall that by Lemma \rf{s2.2}
$$
\eqalign{
    \int_{\GM u} \phi \cdot \nus{\GM u} \, d\H^n
= & \int_\O \big[ \phi^x(x,u) \cdot \nabla u - \phi^t(x,u) \big] \, dx \cr
  &  + \int_{Su} \Big[ \int_{u^-}^{u^+} \phi^x(x,t) \, dt \Big] \cdot \nu_u
     \, d\H^{n-1} \ . \cr
}
\eqno (\frm{3.2})
$$
It is an elementary fact that for every $\xi$, $\eta\in\R^n$
we have $\xi\cdot \eta - {1\over4} |\xi|^2 \le |\eta|^2$, and
equality holds if and only if $\xi=2\eta$.
Hence, setting $\xi:=\phi^x(x,u)$ and $\eta:=\nabla u$,
and taking (a) into account, we obtain that, $\L^n$-a.e.\ on $\O$,
$$
\eqalign{
\phi^x(x,u) \cdot \nabla u-\phi^t(x,u)
& \le \phi^x(x,u) \cdot \nabla u- {1\over4} |\phi^x(x,u)|^2 +
  \beta (u-g)^2 \cr
& \le |\nabla u|^2 + \beta(u-g)^2 \cr
}
$$
and consequently
$$
\int_\O \big[ \phi^x(x,u) \cdot\nabla u -\phi^t(x,u) \big] \, dx
\le \int_\O \big[ |\nabla u|^2 + \beta(u-g)^2 \big]\, dx \ .
\eqno (\frm{3.4})
$$
Moreover, equality holds in (\rf{3.4}) if and only if
$\phi^x(x,u)=2\nabla u$ and
$\phi^t(x,u)={1\over4}|\phi^x(x,u)|^2 - \beta (u-g)^2 =
|\nabla u|^2 - \beta (u-g)^2 $
for a.e.\ $x\in\O$, which is (a').

As for the second integral in the right-hand side of (\rf{3.2}),
condition (b) above implies
$$
\Big[ \int_{u^-}^{u^+} \phi^x(x,t) \, dt \Big] \cdot \nu_u
\le \Big| \int_{u^-}^{u^+} \phi^x(x,t) \, dt \Big| \le \alpha
\quad\hbox{$\H^{n-1}$-a.e.\ on $Su$,}
$$
and then
$$
\int_{Su} \Big[ \int_{u^-}^{u^+} \phi^x(x,t) \, dt \Big]
\cdot \nu_u \, d\H^{n-1} \le \alpha\H^{n-1}(Su) \ .
\eqno (\frm{3.3})
$$
Moreover it is clear that equality holds in (\rf{3.3}) if and only if
(b') is satisfied.

Inequality (\rf{3.1.0}) follows now from (\rf{3.2}), (\rf{3.4}),
and (\rf{3.3}), as well as the rest of the statement.
\qed

\Rem{\thm{s3.1.2}}
Let $\G$ be the class of all vectorfields $\phi$ on $\O\times\R$,
not necessarily bounded, which satisfy assumptions (a) and (b) of Lemma
\rf{s3.1.0} with $U:=\O\times\R$.
It can be easily proved that for every $u$ in $SBV(\O)$ there exists 
$\phi\in\G$ which satisfies
assumptions (a') and (b') for $u$, so that equality holds in (\rf{3.1.0}).
Starting from this, one can also find vectorfields $\phi\in\G$ of class
$C^1_c$ such that the value of the flux in the right-hand side
of (\rf{3.1.0}) is arbitrarily close to $F(u)$. Taking (\rf{2.3})
into account, we can thus prove that for all $u\in SBV(\O)$ there holds 
(cf.\ (\rf{1.6}))
$$
F(u) 
=\sup_{\phi\in\G\cap C^1_c} \int_{\GM u} \phi \cdot \nus{\GM u} \, d\H^n
=\sup_{\phi\in\G\cap C^1_c} \int_{\O\times\R} \phi \cdot D1_u \ .
\eqno (\frm{3.4.1})
$$
Moreover one can show that for any function $u$ which is in $BV(\O)$, 
but not in $SBV(\O)$, the last two terms in (\rf{3.4.1}) are equal to
$+\infty$. Since every integral of the form 
$\int_{\O\times\R} \phi \cdot D1_u$, with $\phi$ of class $C^1_c$ on
$\O\times\R$, is continuous with respect to the weak* topology of 
$BV(\O)$, formula (\rf{3.4.1}) shows that the Mumford-Shah functional 
$F$, extended to $+\infty$ to the rest of
$BV(\O)$, is weak* lower semicontinuous.
The same argument can be easily applied to the general functionals
considered in the next subsection, providing another proof
of the well-known compactness and semicontinuity results in $SBV$ 
due to L. Ambrosio (see [\rf{Amb-ARMA}], or [\rf{AFP}], Sections 4.1 
and 5.4).

\Theor{\thm{s3.2.0}}
{\it
Let $u$ be a function with complete graph contained in $\ove U$,
and assume that there exists an approximately regular vectorfield
$\phi$ on $\ove U$ which is  {\rm divergence-free} on $U$ and satisfies
assumptions {\rm(a)}, {\rm(b)}, {\rm(a')}, and {\rm(b')} of Lemma
\rf{s3.1.0}. Then $u$ is a Dirichlet $\ove U$-minimizer of $F$.
If, in addition, the normal component of $\phi$ at the boundary of
$\O\times\R$ vanishes, namely
$$
\phi^x \cdot \nus{\bd\O}=0
\quad\hbox{$\H^n$-a.e.\ on $(\bd\O\times\R) \cap\bd U$,}
\eqno (\frm{3.4.2})
$$
then $u$ is also an absolute $\ove U$-minimizer of $F$.
}

\Pr
Let $v$ be a function in $SBV(\O)$ such that $v=u$ on
$\bd\O$ and $\GM v\subset\ove U$. Then
$$
F(v)
\ge \int_{\GM v} \phi \cdot \nus{\GM v} \, d\H^n
     = \int_{\GM u} \phi \cdot \nus{\GM u} \, d\H^n = F(u) \ .
\eqno (\frm{3.5})
$$
Here, the first inequality and the last equality follow from
Lemma \rf{s3.1.0}, while the first equality follows from Lemma
\rf{s2.6}. 
We have thus proved that $u$ is a Dirichlet $\ove U$-minimizer of $F$.
Viceversa, assuming (\rf{3.4.2}) we obtain that the first equality
in (\rf{3.5}) holds even if the traces of $v$ and $u$ on $\bd\O$ differ, 
which proves that $u$ is an absolute $\ove U$-minimizer of $F$.
\qed

\Def{\thm{s3.2.1}}
{\it
We call the vectorfield $\phi$ in the first part of Theorem \rf{s3.2.0}
a {\rm Dirichlet calibration for $u$ on $\ove U$} (with respect to $F$).
If $\phi$ satisfies the additional assumption {\rm(\rf{3.4.2})}, then we 
call it an {\rm absolute calibration}.}

\medskip
When $U:=\O\times\R$ we omit to write it.
When it is clear from the context, we may also omit to specify the 
functional, the set $U$, and whether the calibration is Dirichlet or
absolute, and simply say that $\phi$ is a calibration for $u$, or that
$\phi$ calibrates $u$.

\Rem{\thm{s3.5.1}}
If $\phi$ is an absolute calibration for $u$, then it is also an absolute 
calibration for every other minimizer. Indeed, if $F(v)=F(u)$, the first 
inequality in (\rf{3.5}) must be an equality, and by Lemma \rf{s3.1.0} 
this means that $\phi$ satisfies assumptions (a') and (b') for $v$ too.
Similarly, if $\phi$ is a Dirichlet calibration for $u$, then it is 
also a Dirichlet calibration for any other Dirichlet minimizer with
the same boundary values as $u$.

This fact can be sometimes used to prove that the minimizer is unique.
For instance, if $\phi$ calibrates a function $u$ with a negligible
singular set (i.e., $\H^{n-1}(Su)=0$), and the inequality in assumption (b)
is always strict, then we deduce that assumption (b') can only be satisfied
by functions with negligible singular sets, and therefore all minimizers
should have this property. But on this class the functional $F$
is strictly convex (for $\beta>0$, and even for $\beta=0$ in case of
Dirichlet minimizers), and therefore the minimizer must be unique
(see Remarks \rf{e6.1}, \rf{e9.1}, \rf{e17.1}, and Paragraphs \rf{e18} and 
\rf{e19}).

\Rem{\thm{s3.5.2}}
The functional $F_0$ in (\rf{1.2}) is obtained by setting $\alpha:=1$
and $\beta:=0$
in the definition of $F$. In this specific case,
assumptions (a), (b), (a'), and (b') in Lemma \rf{s3.1.0} become

\smallskip
\itemm{(a)}
{$\phi^t(x,t)\ge {1\over 4} |\phi^x(x,t)|^2$
for $\L^n$-a.e.\ $x\in\O$ and every $t\in[\tau_1,\tau_2]$,}

\itemm{(b)}
{$\dps\Big| \int_{t_1}^{t_2} \phi^x(x,t) \, dt \Big|\le 1$
for $\H^{n-1}$-a.e.\ $x\in\O$ and every
$t_1,t_2\in [\tau_1,\tau_2]$,}

\itemm{(a')}
{$\phi^x(x,u)=2\nabla u$ and $\phi^t(x,u)=|\nabla u|^2$ for
$\L^n$-a.e.\ $x\in\O$,}

\itemm{(b')}
{$\dps\int_{u^-}^{u^+} \phi^x(x,t) \, dt=\nu_u$ 
for $\H^{n-1}$-a.e.\ $x\in Su$.}

\Rem{\thm{s3.5.3}}
It must be noticed that, given a boundary value $w$ in 
$L^1(\bd\O,\H^{n-1})$, the Dirichlet problem 
$$
\min\Big\{
\hbox{$F_0(u)$ : $u\in SBV(\O)$, $u=w$ $\H^{n-1}$-a.e.\ on $\bd\O$,
$\GM u\subset\overline U$}
\Big\}
\eqno (\frm{3.7.1})
$$
may not have a solution, even for a very regular $w$, due to a lack of 
continuity of the trace operator on $SBV(\O)$. Therefore problem 
(\rf{3.7.1}) is usually replaced by the relaxed problem 
$$
\min
\Big\{
\hbox{$F_0(u) + \H^{n-1}\big(\{x\in\bd\O: \, u(x)\ne w(x) \}\big)$ 
: $u\in SBV(\O)$, $\GM u\subset\overline U$} \Big\}
\eqno (\frm{3.7.2})
$$
(where in the second term $u$ denotes, as usual, the trace of $u$ on 
$\bd\O$).
A variant of the standard lower semicontinuity and compactness theorems 
in $SBV(\O)$ (see [\rf{Cel}] or [\rf{BCDM}]) shows that problem
(\rf{3.7.2}) has always a solution. 

The calibration method applies to problem (\rf{3.7.2}), too.
In this case calibrations are approximately regular vectorfields 
$\phi$ on $\ove\O\times\R$ which are divergence-free, and
satisfy conditions (a), (b), (a'), (b') of Remark \rf{s3.5.2} 
and,  in addition, the following two conditions
(which may be viewed as extensions of (b) and (b') to the boundary):

\itemm{(c)}
{$\dps\Big| \int_w^s \phi^x(x,t) \, dt \Big|\le 1$
for $\H^{n-1}$-a.e.\ $x\in\bd\O$ and every
$s\in [\tau_1,\tau_2]$,}

\itemm{(c')}
{$\dps\int_w^u \phi^x(x,t) \, dt=\nus{\bd\O}$ 
for $\H^{n-1}$-a.e.\ $x\in\bd\O$ with $u(x)\ne w(x)$,}

\smallskip\noindent
where $u$ and $w$ are computed at $x$, as well as $\tau_1$ and 
$\tau_2$. Indeed, arguing as in Lemma \rf{s3.1.0} we obtain
$$
\eqalign{
F_0(v) +
& \H^{n-1}\big(\{x\in\bd\O: \, v(x)\ne w(x) \}\big) \ge \cr
& \ge\int_{\GM v} \phi \cdot \nus{\GM v} \, d\H^n 
  + \int_{\bd\O} \Big[ \int_w^v \phi^x(x,t) \, dt\Big] \cdot \nus{\bd\O} 
  \, d\H^{n-1} \cr
}
\eqno (\frm{3.7.3})
$$
for every $v\in SBV(\O)$, with equality for $v=u$, moreover the right 
hand side of (\rf{3.7.3}) does not depend on $v$ (apply Lemma \rf{s2.6}
taking into account that $\phi$ is divergence free).

\subsection{Calibrations for general functionals}
The method of calibrations can
be easily adapted to a larger class of functionals
on $SBV(\O)$. Take indeed
$$
\Psi(u):=\int_\O f(x,u,\nabla u)\,dx
+\int_{Su} \psi(x,u^-,u^+,\nu_u) \, d\H^{n-1} \ ,
\eqno (\frm{3.8})
$$
where $f: \O\times\R\times\R^n\to[0,+\infty]$ and
$\psi: \O\times\R\times\R\times \SS^{n-1}\to[0,+\infty]$.
We refer
to [\rf{Amb-ARMA}] for general conditions
on $f$ and $\psi$ which imply the lower semicontinuity of
the functional (\rf{3.8}) and guarantee the existence
of minimizers.
However, lower semicontinuity
is irrelevant for the calibration method.

Let $f^*$ and $\bd_\xi f$ denote the convex
conjugate and the subdifferential of $f$ with respect to
the last variable. We recall that the subdifferential of a
function $g: \R^n\to[0,+\infty]$ at the point $\xi \in\R^n$
is defined as the set of vectors $\eta \in\R^n$ such that 
$g(\xi)+\eta\cdot (\zeta-\xi)\le g(\zeta)$ for every $\zeta\in\R^n$. 
It is well known that for every $\xi$, $\eta\in\R^n$ 
we have the inequality $\xi \cdot \eta -g^{*}(\eta)\le g(\xi)$,
and that equality holds if and only if $\eta \in \bd_\xi g(\xi)$.
Using these properties we obtain the following variant of Lemma 
\rf{s3.1.0}, whose proof is omitted.

\Lemma{\thm{s03.1}}
{\it
Let $\phi$ be a vectorfield on $\ove U$ which
satisfies the following assumptions:

\smallskip
\itemm{(a)}
{$\phi^t(x,t) \ge f^*(x,t,\phi^x(x,t))$
for $\L^n$-a.e.\ $x\in\O$ and every $t\in [\tau_1,\tau_2]$,}

\itemm{(b)}
{$\dps\Big[ \int_{t_1}^{t_2} \phi^x(x,t) \, dt \Big] 
\cdot\nu \le\psi(x,t_1,t_2,\nu)$
for $\H^{n-1}$-a.e.\ $x\in\O$, every $\nu\in \SS^{n-1}$, and every
$t_1<t_2$ in $[\tau_1,\tau_2]$.}

\smallskip\noindent
Then for every $u$ with complete graph contained in $\ove U$ we have
$$
\Psi(u)\ge \int_{\GM u} \phi \cdot \nus{\GM u} \, d\H^n \ .
\eqno(\frm{ab10})
$$
Moreover, equality  holds in {\rm (\rf{ab10})} for a given $u$ 
if and only if

\smallskip
\itemm{(a')}
{$\phi^x(x,u)\in\bd_\xi f(x,u,\nabla u)$ and
$\phi^t(x,u)=f^*(x,u,\phi^x(x,u))$ for $\L^n$-a.e.\ $x\in\O$,}

\itemm{(b')}
{$\dps \Big[ \int_{u^-}^{u^+} \phi^x(x,t) \, dt \Big] \cdot
\nu_u = \psi ( x,u^-,u^+,\nu_u )$ 
for $\H^{n-1}$-a.e.\ $x\in Su$,}

\smallskip\noindent
where $u$, $u^\pm$, $\nabla u$, and $\nu_u$
are always computed at $x$.
}

\medskip
Proceeding as in the previous subsection, one can prove the following
analogue of Theorem \rf{s3.2.0}.

\Theor{\thm{s03.15}}
{\it
Let $u$ be a function with graph contained in $\ove U$.
Assume that there exists an approximately regular vectorfield $\phi$
on $\ove U$ which is {\rm divergence-free} and satisfies
assumptions {\rm(a)}, {\rm(b)}, {\rm(a')}, and {\rm(b')} of Lemma
\rf{s03.1}. Then $u$ is a Dirichlet $\ove U$-minimizer of $\Psi$.
If, in addition, the normal component of $\phi$ on the boundary of
$\O\times\R$ vanishes, i.e., {\rm(\rf{3.4.2})} holds,
then $u$ is also an absolute $\ove U$-minimizer of~$\Psi$.
}

\subsection{Calibrations for minimal partitions}
Besides $F$ and $F_0$, interesting examples of functionals
of the form (\rf{3.8}) arise in
different variants of the minimal partition problem.
Let us consider the case where the number $m$ of the elements of the
partition is prescribed. To formulate the problem, we fix
real numbers $a_1<\cdots <a_m$ in an arbitrary way, and consider
only functions $u$ in the class $BV(\O,\{a_i\})$ of all
$u\in BV(\O)$ which take only these prescribed values.
The corresponding level sets $A_i:=\{u=a_i\}$, sometimes called 
{\it phases\/},
form a partition $\O$; for $i<j$, the set $S_{ij}$ of all $x\in Su$
such that $u^-(x)=a_i$ and $u^+(x)=a_j$ is called the {\it interface\/} 
between the phases $A_i$ and $A_j$, and is oriented by the normal
$\nu_{ij}$ pointing from $A_i$ to $A_j$ (hence $\nu_{ij}=\nu_u$).

We consider functionals of the form
$$
\F(A_1,\ldots,A_m)=\sum_{i<j} \int_{S_{ij}} \hskip -.1cm
\psi_{ij}(x,\nu_{ij}) \, d\H^{n-1}\ ,
\eqno (\frm{3.9})
$$
where $\psi_{ij}:\O\times \SS^{n-1}\to[0,+\infty]$.
Notice that the weights $\psi_{ij}$ may depend on the point at which
two phases meet and on the normal to the interface at that point.
A partition $(A_{1},\ldots,A_{m})$ is said to be a {\it Dirichlet
minimizer\/} of $\F$ if it minimizes $\F$ among all partitions
$(B_{1},\ldots,B_{m})$ such that, for every $i$, the characteristic
functions of $A_{i}$ and $B_{i}$ have the same trace on $\bd\O$.
Notice that all these notions do not depend on the particular
choice of the numbers $a_{1},\ldots,a_{m}$.

We now define a functional $\Psi$ of type (\rf{3.8}) by setting
$$
\eqalign{
f(x,t,\xi)
& :=\cases{0 & if $t\in\{a_i\}$ for some $i$ and $\xi=0$, \cr
     +\infty & otherwise, \cr} \cr
\psi(x,t_1,t_2,\nu)
& :=\cases{\psi_{ij}(x,\nu) & if $t_1=a_i$ and $t_2=a_j$ for some $i<j$,\cr
                   +\infty		& otherwise. \cr} \cr
}
$$
Since every function $u$ in $BV(\O,\{a_i\})$ belongs to $SBV(\O)$, 
and $\nabla u=0$ $\L^n$-a.e.\ in $\O$, one easily checks that
$\Psi(u)$ is finite only if $u$ belongs to $BV(\O,\{a_i\})$, 
and in this case $\Psi(u)=\F(A_{1},\ldots,A_m)$, where $(A_1,\ldots,A_m)$ 
is the partition associated with $u$. Hence a partition $(A_1,\ldots,A_m)$
is a Dirichlet minimizer for $\F$ if and only if the corresponding function 
$u$ is a Dirichlet minimizer for $\Psi$, which according to Theorem 
\rf{s03.15} is implied by the existence of a calibration.

For this particular choice of $\Psi$, a calibration for 
$u\in BV(\O,\{a_i\})$ 
is an approximately regular vectorfield $\phi$ on $\ove\O\times\R$
which is divergence-free and satisfies the following properties 
(cf.\ Lemma \rf{s03.1}):

\smallskip
\itemm{(a)}
{$\phi^t(x,a_i) \ge 0$ for $\L^n$-a.e.\ $x\in\O$ and every $i$,}

\itemm{(b)}
{$\dps\Big[ \int_{a_i}^{a_j} \phi^x(x,t) \, dt \Big]\cdot\nu
\le\psi_{ij}(x,\nu)$
for $\H^{n-1}$-a.e.\ $x\in\O$ and every $\nu\in \SS^{n-1}$ and $i<j$,}

\itemm{(a')}
{$\phi^t(x,a_i)=0$ for $\L^n$-a.e.\ $x\in A_i$ and every $i$,}

\itemm{(b')}
{$\dps \Big[\int_{a_i}^{a_j} \phi^x(x,t) \, dt
\Big]\cdot\nu_{ij}= \psi _{ij}(x,\nu_{ij})$
$\H^{n-1}$-a.e.\ on $S_{ij}$ for every $i<j$.}

\medskip
These calibrations can be re-written in a
different and more interesting way:
let us set, for $i=1,\ldots,m$, $x\in\ove\O$,
$$
\phi_i(x):=\int_{a_1}^{a_i} \hskip-.2cm \phi^x(x,s) \, ds \ .
$$
One can check that the vectorfields $\phi_i$ are approximately 
regular on $\ove\O$ and have divergence in $L^\infty(\O)$
(more precisely, $\div\phi_i(x)=\phi^x(x,a_1)-\phi^x(x,a_i)$
-- apply, e.g., formula (\rf{2.2.3}) with $\O$ replaced
by $\O\times (a_1,a_i)$, and $u$ any smooth function on 
$\O\times (a_1,a_i)$ which depends only on $x$)
and satisfy the following properties:

\smallskip
\itemm{(c)}
{$\div\phi_i\ge\div\phi_j$ for $\L^n$-a.e.\ $x\in A_i$
for every $j\ne i$,}

\itemm{(d)}
{$\big(\phi_j-\phi_i\big)\cdot\nu\le\psi_{ij}(x,\nu)$
for $\H^{n-1}$-a.e.\ $x\in\O$ and every $\nu\in \SS^{n-1}$, $i<j$,}

\itemm{(d')}
{$\big(\phi_j-\phi_i\big)\cdot\nu_{ij} =
\psi_{ij}(x,\nu_{ij})$
for $\H^{n-1}$-a.e.\ $x\in S_{ij}$ and every $i<j$.}

\smallskip\noindent
Conversely, given approximately regular vectorfields
$\phi_i$ on $\ove\O$ with divergence in $L^\infty(\O)$ 
which satisfy (c), (d), (d'),
we define a vectorfield $\phi$ on $\ove\O\times\R$ 
as follows: we take smooth non-negative functions $\sigma_i$
with support included in $(a_i,a_{i+1})$ and integral equal
to $1$, and set
$$
\phi^x(x,t)
:=\sigma_i(t) \,\big(\phi_{i+1}(x)-\phi_i(x) \big)
\quad\hbox{for $x\in\ove \O$, $a_i\le t \le a_{i+1}$,}
$$
then we take $\phi^t$ so that
$$
\cases{
\phi^t(x,a_i):=0 & for $x\in A_i$, \cr
\bd_t\phi^t(x,t):=\sigma_i(t)\, \big( \div\phi_i(x)-\div\phi_{i+1}(x) \big) 
    \vspazio \hskip -1.5 cm \cr
& for $x\in\O$, $a_i\le t \le  a_{i+1}$; \cr
}
\eqno (\frm{3.11})
$$
the definition is completed by setting $\phi(x,t):=\phi(x,a_{1})$ for
$t<a_{1}$ and $\phi(x,t):=\phi(x,a_{m})$ for $t>a_{m}$.

Then $\phi$ is divergence free in $\O\times\R$ by 
construction, and one can easily check 
that it satisfies assumptions (a), (b), (a'), and (b') 
above. Moreover, since each $\phi_i$ is approximately regular on $\ove\O$, 
one can verify that $(\phi^x,0)$ is approximately regular on $\O\times\R$, 
and the same holds for $(0,\phi^t)$ by Remark \rf{r1}. Hence $\phi$ 
is approximately regular too.

We have thus proved the following result.

\Theor{\thm{s03.16}}
{\it
Let $(A_1,\ldots, A_{m})$ be a partition of $\O$, and assume that 
there exist approximately regular vectorfields $\phi_1,\ldots,\phi_m$ 
on $\ove \O$ with divergences in $L^\infty(\O)$ which satisfy assumptions 
{\rm(c)}, {\rm(d)}, and {\rm(d')} above. 
Then $(A_1,\ldots, A_{m})$ is a Dirichlet minimizer
of the functional $\F$ in {\rm (\rf{3.9})}. }

\Rem{\thm{s03.17}}
A particularly relevant example of functional
of type (\rf{3.9}) is the ``interface size'',
which is obtained by taking $\psi_{ij}\equiv 1$ for all $i<j$.
In this case assumptions {\rm(d)} and {\rm(d')} above
reduce to

\smallskip
\itemm{(d)}
{$| \phi_j(x) - \phi_i(x) |\le 1$ for every $x\in\O$ and every $i<j$,}

\itemm{(d')}
{$\phi_j(x) - \phi_i(x)=\nu_{ij}(x)$ for
$\H^{n-1}$-a.e.\ $x\in S_{ij}$ and every $i< j$.}

\smallskip
Calibrations of this type have already been introduced in [\rf{ML1}] and
[\rf{Bra1}] as ``paired calibrations''. More precisely, a paired 
calibration for a partition $(A_1,\ldots,A_m)$ is an ordered $m$-uple of
approximately regular vectorfields $\phi_1,\ldots,\phi_m$ on $\ove\O$
which are divergence-free and satisfy
assumptions (d) and (d') above. We notice that the assumption that
the vectorfields $\phi_i$ are divergence-free is stronger than (c),
and indeed our definition allows in principle  for a larger class of
calibrations.

Among other applications, in [\rf{ML1}] it is shown that in any dimension
$n$ the partition of a regular simplex in $\R^n$ given by the $n+1$ 
simplices spanned by one face and the centre is a Dirichlet minimizer 
of the interface size (and the paired calibration consists simply of $n+1$ 
{\it constant} vectorfields which are orthogonal to the corresponding 
faces).
For $n=3$ this statement was first shown in [\rf{Taylor}] with a 
(relatively) long proof.
In [\rf{Bra1}] it is shown that, unlike what happens in dimension
$3$, the partition of a hypercube in $\R^n$, $n\ge 4$, given by the
$2n$ simplices spanned by one face and the centre is a Dirichlet minimizer
of the interface size.
In both papers the theory is also extended to cover more general 
functionals.
For further examples and results, see also [\rf{Bra2}],
[\rf{Bra3}]; other references are included in the survey~[\rf{Mor}].

\Section{\chp{s4}.~Applications to the homogeneous
\vskip 4 pt
Mumford-Shah functional}

In this section we give some examples of Dirichlet minimizers
of the homogeneous Mumford-Shah functional $F_0$.
We begin with a few remarks which may
be useful when constructing calibrations.

\Rem{\thm{e1}}
By a simple truncation argument, to prove that a function $u: \O\to[m,M]$
is a (Dirichlet) minimizer for $F_0$ it suffices to show that
$F_0(u)\le F_0(v)$ for all competitors $v$ such that $m\le v \le M$.
Thus it is enough to show that $u$ is a Dirichlet $\ove U$-minimizer, with
$U:=\O\times(m,M)$.
In the following we often tacitly assume this principle, and construct
calibrations in $\ove\O\times[m,M]$ instead of $\ove\O\times\R$.
Notice, however, that a calibration $\phi$ on $\ove\O\times [m,M]$ 
can be extended to $\ove\O\times\R$ in a rather simple way:
it suffices to set $\phi(x,t):=\big(0,\phi^t(x,m)\big)$ 
for $t< m$ and $\phi(x,t):=\big(0,\phi^t(x,M)\big)$ for $t> M$ 
(cf.\ Remarks \rf{r1} and \rf{s2.9}(b)).

The same conclusion holds for $F$ if $g$ satisfies $m\le g\le M$ too
(but may fail for a functional of the general form (\rf{3.8}), due 
to lack of suitable truncations).

\Rem{\thm{e4}}
We can construct divergence-free vectorfields on an open set
$A\subset \O\times\R$ using fibrations of $A$ by graphs of harmonic
functions.
More precisely, given harmonic functions $u_\ll$ whose graphs are
pairwise disjoint and cover $A$, for all $(x,t)\in A$ we set
$$
\phi(x,t):= \big( 2\nabla u_\ll (x), |\nabla u_\ll (x)|^2 \big) \ ,
\eqno(\frm{4.2})
$$
where $\ll=\ll(x,t)$ is taken so that $t=u_\ll(x)$.
Thus $\phi$ is a vectorfield on $A$ which, by construction, satisfies
assumption (a) of Remark \rf{s3.5.2}, and assumption (a')
for every $u_\ll$.

We prove that $\phi$ is divergence-free under the additional
assumption that the function $u(x,\ll):=u_\ll(x)$ is of class $C^1$
and $\bd_\ll u(x,\ll)\neq 0$ for every $(x,\ll)$, which implies that 
the parameter $\ll$ can be (locally) chosen so that it depends on $x$ 
and $t$ in a $C^1$ fashion. Then we get
$$
\eqalignno{
\div\phi
& = 2\Delta_x u + 2\bd_\ll \nabla\!_x u \cdot \nabla\!_x \ll
    + 2\nabla\!_x u \cdot \bd_\ll \nabla\!_x u \, \bd_t \ll \cr
& = 2\bd_\ll \nabla\!_x u \cdot (\nabla\!_x \ll + \nabla\!_x u
    \, \bd_t \ll) \ . & (\frm{4.2}) \cr
}
$$
On the other hand, deriving the identity $u(x,\ll(x,t))=t$ with respect
to $x$ and $t$ we get $\nabla\!_x u + \bd_\ll u \, \nabla\!_x\ll=0$
and $\bd_\ll u \, \bd_t\ll=1$, respectively. This implies that the last
factor in (\rf{4.2}) vanishes, and thus $\phi$ is divergence-free
(to make this argument work, we need that $\nabla\!_x u$ is of class $C^1$ 
in $(x,\ll)$, which can be derived by the fact that each function $u_\ll$ 
is harmonic).

In Paragraphs \rf{e6} and \rf{e9} below we apply this idea
by embedding a harmonic function that we intend to calibrate
into a family of harmonic functions whose graphs fibrate
$A:=\O\times(m,M)$, and taking $\phi$ as in (\rf{4.2}). 
To show that $\phi$ is a calibration we will have only to 
verify assumption (b) of Remark \rf{s3.5.2}.

\Rem{\thm{e4.1}}
The construction described in the previous remark is a particular
case of a classical result about extremal fields of
scalar functionals (see, e.g., [\rf{AAC}], Section 4).
Let be given an open subset $A$ of $\O\times\R$ which is covered
by a family of pairwise disjoint graphs of solutions $u_\ll$ of
the Euler-Lagrange equation associated with the functional
$\int_\O f(x,u,\nabla u) \, dx$, that is
$$
\div\big(\bd_\xi f(x,u,\nabla u)\big)=\bd_t f(x,u,\nabla u) \ .
$$
Here $t$ and $\xi$ denote the second and third variable
in the argument of $f(x,t,\xi)$, corresponding to $u$ and $\nabla u$,
and we assume that $f$ is of class $C^2$ in $(t,\xi)$
and convex in $\xi$.
For all $(x,t)\in A$ let
$$
\cases{
\phi^x(x,t):= \bd_\xi f(x,u_\ll(x),\nabla u_\ll(x)) \ ,\cr
\phi^t(x,t):= \bd_\xi f(x,u_\ll(x),\nabla u_\ll(x)) \cdot \nabla u_\ll(x)
-f(x,u_\ll(x),\nabla u_\ll(x)) \vspazio\cr
}
\eqno (\frm{4.2.1})
$$
where $\ll=\ll(x,t)$ is taken so that $t=u_\ll(x)$.
Then $\phi$ satisfies by construction assumption (a) 
of Lemma \rf{s03.1}, and assumption (a') for every function $u_\ll$.
Moreover one can prove that $\phi$ is divergence-free
under the additional assumption that the
function $u(x,\ll):=u_\ll(x)$ is of class $C^2$ in $(x,\ll)$ and 
$\bd_\ll u\neq 0$ (cf.\ [\rf{AAC}], Theorem 4.6).
If, in addition, $\phi$ satisfies assumption (b) of
Lemma \rf{s03.1}, then it is a calibration for each
$u_\lambda$ relative to the functional $\Psi$ in (\rf{3.8}), 
and hence each $u_\ll$ is a $\ove U$-Dirichlet minimizer of $\Psi$ 
provided that $\ove U$ is contained in $A$.

\Rem{\thm{e4.2}}
The construction described in the previous
remarks cannot be really used for absolute minimizers, i.e.,
when no Dirichlet boundary condition is imposed.
Indeed, calibrations for absolute minimizers should have
vanishing normal component at the boundary of $\O\times\R$
(see Theorem \rf{s03.15}), and this holds true for the
vectorfield $\phi$ in (\rf{4.2.1}) if and only if the functions
$u_\ll$ satisfy the natural boundary condition
$\bd_\xi f(x,u_\ll(x),\nabla u_\ll(x)) \cdot \nus{\bd\O}(x)=0$ for
$x\in \bd\O$.
But this means that all $u_\ll$ solve the Neumann
problem associated with the functional $\int_\O f(x,u,\nabla u) \, dx$,
while in general we cannot expect a one-parameter family of
solutions for such a problem. For instance, if this functional 
is strictly convex (and this is indeed the case for the regular part of 
the Mumford-Shah functional $F$ when $\beta>0$), then the associated 
Neumann problem admits at most one solution.

\Rem{\thm{e5}}
For $n=1$ and $\O=(a,b)$, the equation $\div\phi=0$ on ${\O\times\R}$,
coupled with the identity $\phi^t={1\over 4}(\phi^x)^2$, reduces
to the first order equation
$$
\textstyle \bd_x \phi^x + {1\over 2} \phi^x \, \bd_t \phi^x=0 \ .
\eqno(\frm{an2})
$$
It easily follows from the method of characteristics that $\phi$
is a $C^1$ solution of (\rf{an2}) in ${\O\times\R}$ if and only if every
level set $\{ \phi^x=s \}$ is composed of straight lines with slope $s/2$
(intersected with ${\O\times\R}$).
In other words, for $n=1$ all $C^{1}$ divergence-free vectorfields
$\phi$ on ${\O\times\R}$ which satisfy $\phi^t={1\over 4}(\phi^x)^2$
(cf.\ conditions (a) and (a') in Remark \rf{s3.5.2})
are associated with a fibration of $\O\times\R$ with
graphs of affine --i.e., harmonic-- functions as in Remark \rf{e4}.

\bigskip
For the rest of this section, calibrations are always intended 
as Dirichlet calibrations for $F_0$, in the sense of Remark \rf{s3.5.2}.
We begin with a discussion of some one-dimensional examples.
Of course, in these examples minimality can be easily checked
by direct computations, and there would be no need for calibrations.
Nevertheless, it is instructive to see what happens, and moreover some
one-dimensional constructions are carried over to higher dimensions
(cf.\ Paragraphs \rf{e9} and \rf{e10}).

\Parag{\thm{e6}.~Affine function in one dimension}
Let $\O$ be the open interval $(0,a)$ and let $u$ be the linear
function $u(x):=\ll x$, with $\ll>0$.
It is easy to see that $u$ is a Dirichlet minimizer of
$F_0$ if and only if
$$
a\ll^{2}\le 1 \ .
\eqno(\frm{4.4.1})
$$
In this case a calibration is given
by the  piecewise constant vectorfield:
$$
\phi(x,t):=\cases{
(2\ll, \ll^{2})
      & if ${\ll\over 2} x\le t \le {\ll\over 2}(x+a)$, \cr
(0,0) \vspazio
      & otherwise.\cr
}
\eqno(\frm{4.5})
$$
Thus $\phi$ satisfies assumptions (a) and (a') of Remark \rf{s3.5.2},
and vanishes outside a stripe of constant height (in grey in Figure 2 below, 
on the left) which is arranged so that (b) holds and $\div\phi$ vanishes 
(cf.\ Remark \rf{s2.9}(b)), while (b') is trivially satisfied.

\figureps{2}{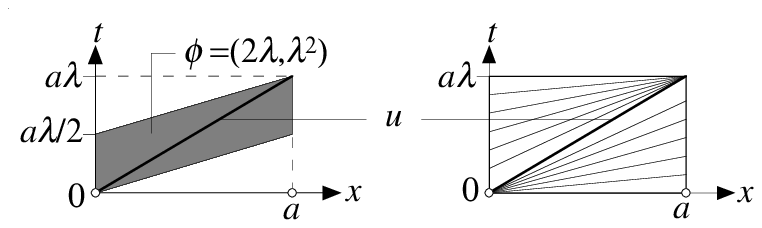}

\noindent
Another calibration is obtained by fibrating the rectangle
$U=(0,a)\times(0,\ll a)$ with affine functions as shown in Figure 2,
on the right, and applying the construction of Remark \rf{e4}:
$$
\phi(x,t):=\cases{
\big( 2{t \over x}, \big({t\over x} \big)^2 \big)
       & if $0\le t\le\ll x$, \cr
\big( 2{\ll a -t \over a-x}, \big( {\ll a -t\over a-x} \big)^2\big) \vspazio
       &if $\ll x\le t\le \ll a$. \cr
}
\eqno(\frm{4.6})
$$
It remains to check that assumption (b) is satisfied, which
happens if and only if (\rf{4.4.1}) holds.

\Rem{\thm{e6.1}}
If $a\ll^2<1$, then both calibrations described in the previous 
paragraph satisfy the strict inequality in assumption (b) of Remark 
\rf{s3.5.2}, i.e.,
$$
\Big| \int_{t_1}^{t_2} \phi^x(x,t) \, dt \Big| < 1
\quad\hbox{for every $x\in [0,a]$ and every $t_1, t_2\in\R$.}
$$
By Remark \rf{s3.5.1}, this shows that the function
$u(x):=\lambda x$ is the unique Dirichlet minimizer of $F_0$ with
$u(0)=0$ and $u(a)=\ll a$.

\Parag{\thm{e7}.~Step function in one dimension}
In Paragraph \rf{e6}, in the limit case $a\ll^2=1$ the linear
function $u(x)=\ll x$ and any step function of the form
$u(x):=0$ for $0<x<c$ and $u(x):=\ll a=\sqrt a$ for $c<x<a$
(with $0<c<a$) are both Dirichlet minimizers with the same boundary 
values. Hence both vectorfields (\rf{4.5}) and (\rf{4.6}) calibrate
these step functions when $\ll:=1/ \sqrt a$ (cf.\ Remark \rf{s3.5.1}).
Furthermore, it is easy to check
that they also calibrate any step function $u$ given by
$u(x):=0$ for $0<x<c$, and $u(x):=h$ for
$c<x<a$ with $h\ge\sqrt a$.

\Rem{\thm{e8}}
When $a\ll^2>1$ the linear function $u(x):=\ll x$ is not a Dirichlet 
minimizer of $F_0$ (a step function is preferable), but it is still a 
Dirichlet $\ove U$-minimizer, when $U$ is the stripe of all points
$(x,t)$ between the graph of $\ll x - {1\over 4\ll}$ and
$\ll x + {1\over 4\ll}$.
A calibration is given by $\phi(x,t):=(2\ll, \ll^{2})$.

\figureps{3}{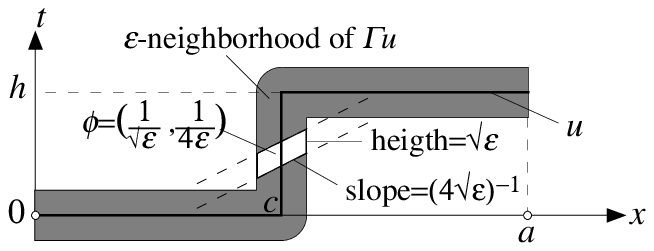}

Conversely, when $h<\sqrt a$, the step function $u$ in Paragraph \rf{e7} 
is no longer a Dirichlet minimizer, but it is Dirichlet $\ove U$-minimizer
when $U$ is an $\eps$-neighbourhood of the complete graph of $u$
(in grey in Figure 3) and $\eps$ satisfies 
${3\over 2} \sqrt\eps+ 2\eps \le h$.
A calibration is given by the piecewise constant vectorfield
which vanishes outside the white parallelogram in Figure 3, 
and is equal to $\big( {1\over\sqrt\eps},{1\over 4\eps} \big)$ inside.

\Parag{\thm{e9}.~Harmonic functions in dimension $n$}
Let $u$ be a harmonic function on $\O$. Since $u$ is
a Dirichlet minimizer of $\int_\O |\nabla u|^2$, it is natural
to ask when it is also a Dirichlet minimizer of $F_0$.
As pointed out by A. Chambolle, this happens when
$$
\osc{u} \cdot \sup |\nabla u| \le 1 \ ,
\eqno (\frm{4.6.1})
$$
where $\osc{u}$ is the oscillation of $u$, namely the difference 
between the supremum $M$ and infimum $m$ of $u$ (over $\O$).
In the one-dimensional case $n=1$ this condition reduces
to (\rf{4.4.1}).

A calibration can be constructed by analogy with (\rf{4.5});
see Figure 4, on the left:
$$
\phi(x,t):=\cases{
\big( 2\nabla u(x), |\nabla u(x)|^2 \big)
			    & if ${1\over 2}(u(x)+m) \le t \le {1\over 2} (u(x)+M)$, \cr
\vphantom{\vrule width 0pt height 13pt depth 0pt} (0,0)
         & otherwise. \cr
}
\eqno (\frm{4.7})
$$
Another calibration can be obtained, as the one in (\rf{4.6}),
by embedding $u$ in a family of harmonic functions whose graphs
fibrate the cylinder $\O\times[m,M]$. More precisely we take
the functions $m+\ll(u-m)$ and $M+\ll(u-M)$ with
$\ll$ ranging in $[0,1]$ (see Figure 4, on the right), and
then the  construction of Remark \rf{e4} gives
$$
\phi(x,t):=\cases{
\big(2{t-m\over u(x)-m}\nabla u(x), ({t-m\over u(x)-m})^2 
|\nabla u(x)|^2 \big)
        & if $m\le t\le u(x)$, \cr
\vphantom{\vrule width 0pt height 13pt depth 0pt}
\big(2{M-t\over M-u(x)} \nabla u(x), ({M-t\over M-u(x)})^2
|\nabla u(x)|^2 \big)
        & if $u(x)\le t\le M$. \cr
}
\eqno (\frm{4.8})
$$
One easily checks that both vectorfields are divergence-free
(see Remarks \rf{s2.9}(b) and \rf{e4}),
and that assumptions (a) and (a') of Remark \rf{s3.5.2} are
satisfied, assumption (b') is always trivially
satisfied, while assumption (b) holds if an only if (\rf{4.6.1}) holds.

\figureps{4}{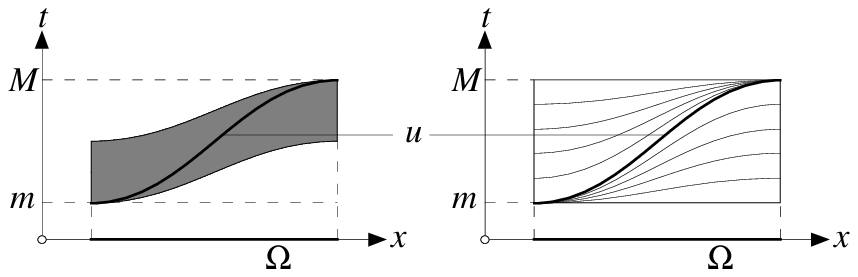}

When (\rf{4.6.1}) is not satisfied, $u$ is still a Dirichlet 
$\ove U$-minimizer of $F_0$, where $U$ is the set of all points 
$(x,t)\in\O\times\R$ which lie between the graph of 
$u(x) - (4|\nabla u(x)|)^{-1}$ and $u(x) + (4|\nabla u(x)|)^{-1}$; 
a calibration is given by
$\phi(x,t):=\big( 2\nabla u(x), |\nabla u(x)|^2 \big)$.

\Rem{\thm{e9.1}}
If inequality (\rf{4.6.1}) holds and $u$ is not affine, 
then the maximum principle implies that $\osc{u} \cdot|\nabla u(x)|<1$ 
for every $x\in\O$, and therefore both calibrations constructed 
in the previous paragraph satisfy the strict inequality in assumption (b) 
of Remark \rf{s3.5.2}.  
By Remark \rf{s3.5.1}, this proves that the harmonic function
$u$ is the only Dirichlet minimizer of $F_0$ with the same boundary 
values as $u$ (cf.\ Remark \rf{e6.1}).

\Parag{\thm{e10}.~Step function in dimension $n$}
Let $\O$ be a product of the form $(0,a)\times V$, where $V$
is a regular domain in $\R^{n-1}$, $n\ge 2$, and let $u$
be the step function given by $u(x):=0$ for $0<x_1<c$,
and $u(x):=h$ for $c<x_1<a$, where $x_1$ denotes the first
coordinate of $x$, $c\in(0,a)$, and $h\ge\sqrt a$.

Using the results in Paragraph \rf{e7} and a symmetrization argument,
it is easy to see that $u$ is a Dirichlet minimizer of $F_0$.
Calibrations can be constructed starting from the one-dimensional
ones described in Paragraph \rf{e7}. More precisely, we take
the vectorfield on $\O\times\R$ which is parallel to the
$(x_1,t)$-plane and is given by formula
(\rf{4.5}) (or even (\rf{4.6})) with $x$ replaced by $x_1$ and
$\ll:=1/\sqrt a$.

\Rem{\thm{e10.1}}
The previous result can be restated by saying that a
step function $u$ with jump of height $h$ along a hyperplane
orthogonal to the direction $e$ is a Dirichlet minimizer whenever
$h\ge\sqrt a$, where $a$ is the diameter of the projection of $\O$
along the $e$ axis.
As in dimension one, any step function
is a Dirichlet $\ove U$-minimizer, when $U$ is an
$\eps$-neighbourhood of the complete graph of $u$ and $\eps$
satisfies ${3\over2}\sqrt{\eps} + 2\eps\le h$ (cf.\ Remark \rf{e8}).
However, unlike what happens in dimension one, $u$ may be a
Dirichlet minimizer of $F_0$ even when the assumption $h\ge\sqrt a$
is not satisfied (cf.\ Paragraph \rf{e15} below).

\Parag{\thm{e11}.~Triple junction in the plane}
Let $\O:=B(0,r)$ be the open disk in the plane with radius
$r$ and centred at the origin, and let $(A_1,A_2,A_3)$ be
the partition of $\O$ defined as follows: $A_i$
is the set of all $x\in\O$ of the form
$x=(\rho\cos\theta, \rho\sin\theta)$,
with ${2\over 3}\pi(i-1)\le\theta<{2\over 3}\pi i$.
Finally define $u:=a_i$ in each $A_i$,
where $a_1 < a_2 < a_3$ are distinct constants.

Thus the singular set of $u$ is given by three line segments
$S_{1,2}$, $S_{2,3}$, and $S_{3,1}$ meeting at the origin with
equal angles (see Figure 5, on the left), 
and it is well-known that this is a minimal network,
in the sense that the corresponding partition $(A_1,A_2,A_3)$
is a Dirichlet minimizer of  the ``interface size'' functional
(see Remark \rf{s03.17}).
Therefore it is natural to conjecture that, when the values
of the constant $a_i$ are sufficiently far apart,
$u$ is a Dirichlet minimizer of $F_0$ too, that is, there is no
convenience in removing part of the jump and taking a function with
non-vanishing gradient.

We prove this conjecture by calibration.
Inspired by the constructions described in Paragraphs \rf{e7} 
and \rf{e10}, we take $e_\pm:=(\pm\sqrt 3/2,-1/2)$, $\ll>0$, 
and set
$$
\phi(x,t):=\cases{
(2\ll e_-, \ll^2)
      & if $|t- {1\over 2}(a_1+a_2)-{\ll\over 2} x\cdot e_- | 
        <{1\over 4\ll}$, \cr
(2\ll e_+, \ll^2) \vspazio
      & if $|t- {1\over 2}(a_2+a_3)-{\ll\over 2} x\cdot e_+ | 
        <{1\over 4\ll}$, \cr
(0,0) \vspazio
      & otherwise. \cr
}
\eqno (\frm{4.9})
$$
Thus $\phi$ is piecewise constant, satisfies assumption (a) of
Remark \rf{s3.5.2} by construction, and vanishes out of two slabs
of constant height $1\over 2\ll$ (see Figure 5, on the right).
These slabs have been arranged in order to fulfill the following
requirements:

\itemm{(i)}
{one slab is contained in $\O\times[a_1,a_2]$ and
the other one in $\O\times[a_2,a_3]$, so that assumption (a') of
Remark \rf{s3.5.2} is satisfied; it is possible to construct such slabs if
we can choose $\ll$ so that $a_{i+1}-a_i \ge\ll r + {1\over 2\ll}$,
that is, if 
$$
a_{i+1}-a_i \ge\sqrt{2r} \ ;
\eqno (\frm{4.9.1})
$$
}

\itemm{(ii)} 
{the compatibility condition (\rf{2.2.2}) is
satisfied on the boundary of the slabs,
so that $\phi$ is approximately regular and
divergence-free (cf.\ Remark \rf{s2.9}(b));
}

\itemm{(iii)} 
{assumption (b') is satisfied for all
points $x$ in $S_{1,2}$ and $S_{2,3}$, where $\nu_u$ coincides with
$e_-$ and $e_+$ respectively.
}

\figureps{5}{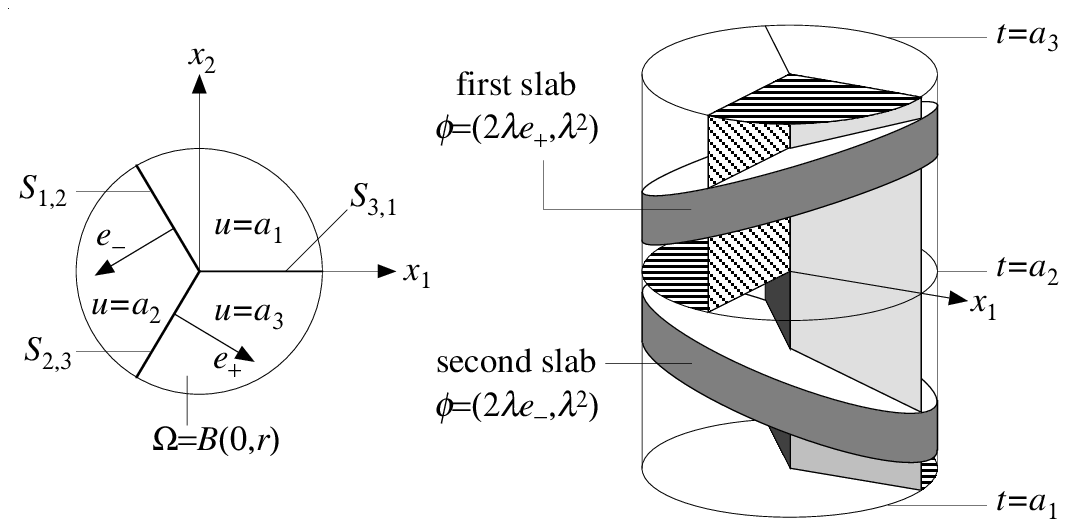}

Moreover (b') holds also for $x$ in $S_{3,1}$, because $e_- + e_+=\nu_u$. 
One also checks that the integral $\int_{t_1}^{t_2} \phi^x(x,t)\,dt$ can 
be always written as a linear combination $\mu_- e_- + \mu_+ e_+$
with $\mu_\pm$ in $[0,1]$ (depending on $x,t_1,t_2$), and
since $e_+$ and $e_-$ span an angle equal to $2\pi/3$, this implies that
the integral has modulus not larger than~1. Thus (b) holds, too.

\Rem{\thm{e12}}
When $a_2-a_1$ (or $a_3-a_2$) is sufficiently small, $u$ is not a 
Dirichlet minimizer. More precisely, if $a_2-a_1 \le{1\over 3}\sqrt r$
(cf. (\rf{4.9.1})), a comparison function $v$ with the same boundary 
values as $u$ and such that $F(v)<F(u)$ is given, in polar coordinates, by
$$
v:=\cases{
{1\over 2}(a_1+a_2)
        & if $0 \le\theta < {4\over 3}\pi$ and $\rho\le r-d$, \cr
{1\over 2}(a_1+a_2) +{1\over 2}(a_1-a_2) \, {1\over d}(\rho-r+d)
        & if $0 \le\theta < {2\over 3}\pi$ and $\rho>r-d$, \vspazio\cr
{1\over 2}(a_1+a_2) +{1\over 2}(a_2-a_1) \, {1\over d}(\rho-r+d)
        & if ${2\over 3}\pi\le\theta<{4\over 3}\pi$ and $\rho>r-d$, 
          \vspazio\cr
a_3
        & if ${4\over 3}\pi \le\theta<2\pi$, \vspazio\cr
}
$$
where $d:=(a_2-a_1)\sqrt r$ (we leave the computations to the reader).

\Parag{\thm{e13}.~Minimal partitions in dimension $n$}
One can generalize the example of the triple junction, and conjecture
the following: if a partition $(A_1,\ldots,A_m)$ of $\O$ is a
Dirichlet minimizer of the ``interface size'' (see Remark \rf{s03.17})
and $u$ is a function which takes a constant value $a_i$ on each $A_i$
(with $a_1<a_2<\ldots<a_m$),
then $u$ is a Dirichlet minimizer of $F_0$ when the values $a_i$ are
sufficiently far apart from each other.
Unfortunately we can only prove this statement under two additional 
assumptions:

\itemm{(i)} 
{the partition $(A_1,\ldots,A_m)$ is not only minimal, but 
admits a paired calibration in the sense of [\rf{ML1}] and [\rf{Bra1}], 
namely there exist approximately regular, divergence-free vectorfields
$\phi_1,\ldots,\phi_m$ on $\ove\O$ which satisfy assumptions (d) and (d') 
in Remark \rf{s03.17};}

\itemm{(ii)}
{for $i=1,\ldots,m-1$ there exist Lipschitz functions
$\psi_i: \ove\O\to\R$ which satisfy almost everywhere the first 
order equation
$$
\textstyle
\nabla \psi_i \cdot (\phi_{i+1}-\phi_i) = {1\over 2}
|\phi_{i+1}-\phi_i|^2 \ .
\eqno (\frm{4.11})
$$
}

\noindent
Adding, if needed, a constant to $\psi_i$, we may also assume that
$$
\osc{\psi_i}=2\|\psi_i\|_infty \ .
\eqno (\frm{g37})
$$
For $i=1,\ldots,m-1$ we take slabs
$U_i$, included in $\ove \O\times(a_i,a_{i+1})$, of the form
$$
\textstyle
U_i=\Big\{ (x,t):
\big| t-{1\over 2}(a_i+a_{i+1}) -\ll_i \psi_i(x) \big|< {1\over 4\ll_i} 
\Big\}\ ,
\eqno (\frm{4.12})
$$
where the constants $\ll_i$ will be chosen below.
Then we set (cf.\ (\rf{4.9}))
$$
\phi(x,t):=\cases{
\big( 2\ll_i(\phi_{i+1}(x)-\phi_i(x)), \ll_i^2 
|\phi_{i+1}(x)-\phi_i(x)|^2 \big) 
   \hskip -3 cm \cr
\vspazio & if $(x,t)\in U_i$ for some $i$, \cr
(0,0) \vspazio & otherwise. \cr
}
$$
Taking into account assumption (d') in Remark \rf{s03.17} and the
definition of the  slabs $U_i$, one can easily check that assumptions 
(a), (a'), and (b') of Remark \rf{s3.5.2} are satisfied. 

Let us check assumption (b). Taken 
$t_1\in[a_i, a_{i+1}]$ and $t_2\in [a_j, a_{j+1}]$ for some $i,j$, 
the integral $\int_{t_1}^{t_2} \phi^x(x,t)\, dt$ can be decomposed 
as the sum of the integrals on the (oriented) intervals 
$[t_1,a_{i+1}]$, $[a_{i+1},a_j]$,
and $[a_j, t_2]$, and hence it can be written as 
$$
\mu_1 (\phi_{i+1}(x)-\phi_i(x)) + (\phi_j(x)-\phi_{i+1}(x)) 
+\mu_2 (\phi_{j+1}(x)-\phi_j(x))
$$
for suitable $\mu_1, \mu_2\in [0,1]$. But this sum can be reorganized as 
the difference between $\mu_2 \phi_{j+1}(x) + (1-\mu_2)\phi_j(x)$ and 
$\mu_1 \phi_i(x) + (1-\mu_1)\phi_{i+1}(x)$. Therefore its modulus 
is the distance between two points in the convex hull
of the vectors $\phi_1(x),\dots, \phi_m(x)$, which has diameter $1$ 
because of assumption (d) in Remark \rf{s03.17}, and (b) is proved.

Since the vectorfields $\phi_i$ are divergence-free and 
approximately regular by assumption, $\phi$ is divergence-free
and approximately regular within each slab (the approximate regularity 
of $(\phi^x,0)$ is immediate, that of $(0,\phi^t)$ follows from 
Remark \rf{r1}),  as well as in the interior of the
complement of the union of all slabs. 
Thus $\phi$ is divergence-free and approximately regular 
in $\O\times\R$ if (and only if) the compatibility condition
(\rf{2.2.2}) is satisfied on the boundary of each slab 
(cf.\ Remark \rf{s2.9}(b)), which reduces to equation (\rf{4.11}).

Therefore we have constructed a calibration for $u$, provided that
we can choose $\ll_i$ so that the slabs $U_i$ are
contained in $\O\times(a_i,a_{i+1})$, that is,
$$
{a_{i+1} - a_i\over 2} 
\ge \ll_i \|\psi_i\|_\infty  + {1\over 4\ll_i}
={\ll_i\over 2} \osc{\psi_i} + {1\over 4\ll_i} \ , 
$$
and this can be done as long as
$$
a_{i+1} - a_i \ge
\sqrt{2\,\osc{\psi_i}}
\quad\hbox{for $i=1,\ldots,m-1$.}
\eqno (\frm{4.13})
$$

\Rem{\thm{e13.1}}
A paired calibration for the partition $(A_1,A_2,A_3)$ described 
in Paragraph \rf{e11} is given by the constant vectorfields
$\phi_1:={1\over 6}(\sqrt 3, 3)$,
$\phi_2:={1\over 6}(-2\sqrt 3,0)$,
$\phi_3:={1\over 6}(\sqrt 3,-3)$, 
and the linear functions $\psi_1$ and $\psi_2$ with
gradients ${1\over 4}(-\sqrt 3,-1)$ and ${1\over 4}(\sqrt 3,-1)$
satisfy equation (\rf{4.11});
if we thus apply the construction of the previous paragraph,
we obtain exactly the calibration described in Paragraph \rf{e11}.

\Rem{\thm{e14}}
The first order equation (\rf{4.11}) does not always admit solutions. 
For instance, since the derivative of $\psi_i$ along the integral 
curves of the vectorfield $\phi_{i+1}-\phi_i$
(i.e., the maximal solutions of the differential equation
$\dot\gamma=\phi_{i+1}(\gamma)-\phi_i(\gamma)$) is always positive,
when there exists a {\it nontrivial\/} closed integral curve within $\O$,
(\rf{4.11}) admits no solution. On the other hand, if $\phi_{i+1}-\phi_i$
is $C^1$ and nowhere vanishing,
and all integral curves start and end at the boundary of $\O$
and intersect a fixed $(n-1)$-dimensional closed manifold $M$ in $\O$
which is transversal to the vectorfield $\phi_{i+1}-\phi_i$, 
i.e., a cross-section of the associated flow, then the method 
of characteristics provides a solution $\psi_i$ to (\rf{4.11}) 
of class $C^1$.

However, such a strong requirement on $\phi_{i+1}-\phi_i$ is far from 
being necessary.
Not only there may exist Lipschitz functions $\psi_i$ which satisfy 
(\rf{4.11}) almost everywhere even if $\phi_{i+1}-\phi_i$ vanishes
somewhere, but  for the purposes of
the previous construction we can even allow $\psi_i$ to be discontinuous
along some integral curve $\gamma$: in this case the boundary of the
slab $U_i$ in (\rf{4.12}) is not just the union of the graphs of
$\ll_i\psi_i+{1\over 2}(a_1+a_2)+{1\over 4\ll_i}$ and
$\ll_i\psi_i+{1\over 2}(a_1+a_2)-{1\over 4\ll_i}$, but there is an 
additional vertical piece contained in $\gamma\times\R$. Yet the
compatibility condition (\rf{2.2.2}) is satisfied there, and then 
$\phi$ is still divergence-free and approximately regular 
(cf.\ Remark \rf{s2.9}(b)).

\Parag{\thm{e15}.~More on the step function in the plane}
Let $\O$ be the rectangle $(-a,a)\times (-b,b)$ in the plane, and
let $(A_1,A_2)$ be the partition of $\O$ given by the sets of all points
$x=(x_1,x_2)$ such that $x_1<0$ and $x_1\ge 0$ respectively.
This partition is obviously minimal,
and a paired calibration is $\phi_1:=(0,0)$, $\phi_2:=(1,0)$.
Setting $p_+:=(0,b)$ and $p_-:=(0,-b)$, another paired calibration 
is given by $\phi_1:=(0,0)$ and
$$
\phi_2(x):=\cases{
(-\sin \theta_+,\cos\theta_+) & for $x\in B(p_+,b)$, \cr
(\sin \theta_-,-\cos\theta_-) \vspazio & for $x\in B(p_-,b)$, \cr
(0,0) \vspazio & otherwise, \cr
}
\eqno (\frm{phi2})
$$
where $\rho_\pm, \theta_\pm$ are the polar coordinates around
the points $p_\pm$; see Figure 6.

\figureps{6}{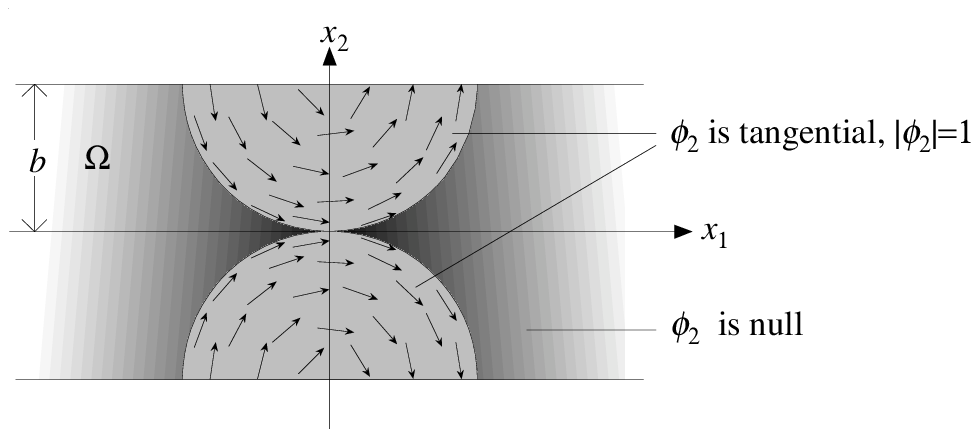}

\noindent
In this case, a function $\psi_1$ which satisfies (\rf{4.11}) almost 
everywhere is given by
$$
\psi_1(x):=\cases{
{1\over 2}(\theta_+ +\pi/2)\rho_+ & for $x\in B(p_+,b)$, \cr
{1\over 2}(\theta_- -\pi/2)\rho_- \vspazio & for $x\in B(p_-,b)$, \cr
0 \vspazio & otherwise, \cr
}
$$
and the construction in Paragraph \rf{e13}, performed with some care
because of the discontinuity of $\psi_1$ along the circles 
$\bd B(p_\pm,b)$, yields a calibration for the step function $u$ which 
takes the value $a_1$ on $A_1$ and $a_2$ on $A_2$, provided that
$a_2-a_1\ge \sqrt{\pi b}$.  Note that the calibration obtained in this
way is defined on the whole stripe $\R\times (-b,b)$ and does not 
depend on~$a$.
This extends the minimality result proved in Paragraph \rf{e10}.

\Parag{\thm{e16}.~More on the triple junction}
Let us apply again the construction described in Paragraph \rf{e13}
to the situation described in Paragraph \rf{e11},
with $\O$ replaced by an $\eps$-neighbourhood of $Su$ within
the ball $B(0,r)$ (in grey in Figure 7).
As already noticed in Remark \rf{e13.1}, a paired calibration
for the partition $(A_2,A_2,A_3)$ is given by the constant
vectorfields $\phi_1:={1\over 6}(\sqrt 3,3)$,
$\phi_2:={1\over 6}(-2\sqrt 3,0)$, and $\phi_3:={1\over 6}(\sqrt 3,-3)$,
but we can take solutions $\psi_i$ of (\rf{4.11}) such that
$|\psi_i| \le 2\eps$ on $\ove\O$, independently of the value of $r$.
More precisely, we take the solution $\psi_i$ of (\rf{4.11})
which takes the value $0$ on the transversal set $M_i$
described in the Figure 7 for $i=1$
(the construction is similar for $i=2$).

\figureps{7}{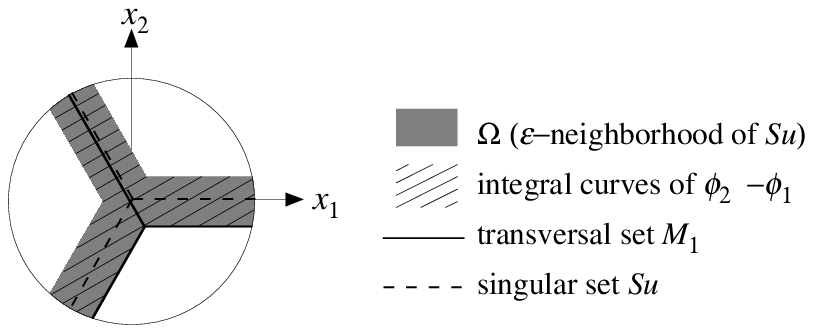}

\noindent
Therefore, if $a_{i+1}-a_i \ge 2\sqrt{2\eps}$,
the construction of Paragraph \rf{e13} yields a calibration of $u$ on
${\ove\O\times\R}$.
Moreover this calibration can be extended to a calibration
on an $\eps$-neighbourhood of the complete graph of $u$ over the
entire $B(0,r)$ by setting it equal to $0$ where it is not yet
defined. Clearly, this requires the slabs $U_i$ are contained in
$B(0,r)\times (a_i+\eps,a_{i+1}-\eps)$ for
$i=1,2$, that is
$$
a_{i+1}-a_i \ge 2\eps + 2\sqrt{2\eps}
\quad\hbox{for $i=1,2$.}
\eqno (\frm{4.15})
$$
This shows that $u$ is a Dirichlet $\ove U$-minimizer of $F_0$
when $U$ is an $\eps$-neighbourhood of the complete graph of
$u$ within $B(0,r)\times\R$ and $\eps$ satisfies (\rf{4.15}).
As expected, $\eps$ does not depend on the size of the domain, 
but only on the relative distances of the values $a_i$.

\Section{\chp{s5}.~Applications to the complete Mumford-Shah
\vskip 4 pt 
functional}
In this section we focus on minimizers of
the complete Mumford-Shah functional $F$ defined in (\rf{1.1})
for $\alpha, \beta>0$, and calibrations will always be intended as
in the sense of Theorem \rf{s3.2.0} and Definition \rf{s3.2.1}.
All the following examples are in dimension $n$.

\Parag{\thm{e17}.~Solutions of the Neumann problem}
If we restrict $F$ to functions of class
$W^{1,2}$, we obtain the strictly convex and coercive functional
$\int_\O [|\nabla u|^2 + \beta(u-g)^2]\, dx$, and its unique minimizer
$u$ is the solution of the Neumann problem
$$
\cases{
   \Delta u=\beta(u-g)  & on $\O$,    \cr
   \bd_\nu u=0 \vspazio & on $\bd\O$, \cr
}
\eqno (\frm{4.16})
$$
where $\bd_\nu$ denotes the normal derivative.
Thus it is natural to ask under which assumptions (on $g$
and $\beta$) $u$ is also a minimizer of $F$ on $SBV(\O)$.
This question is akin to the minimality of harmonic functions 
for $F_0$ discussed in Paragraph \rf{e9}, and following 
the same ideas we can construct an absolute calibration for $u$
provided that $u$ satisfy condition (\rf{4.17.1}) below (cf.\ also
Paragraph \rf{e18}).

More precisely, we assume that $u$ is of class $C^1$ up to the boundary
(this is always satisfied if $\bd\O$ is of class $C^{1,\eps}$ for
some $\eps>0$), we denote the infimum and
the supremum of $g$ by $m$ and $M$ respectively, and set
$$
A:=\Big\{
(x,t)\in\ove\O\times\R: \ {u(x)+m \over 2} \le t\le {u(x)+M \over 2}
\Big\} \ ,
$$
and (cf. (\rf{4.7}))
$$
\phi^x(x,t):=\cases{
    2\nabla u(x) & if $(x,t)\in A$, \cr
    0 \vspazio & otherwise. \cr
}
\eqno (\frm{4.17})
$$
Note that, by the maximum principle, $m\le u\le M$ on $\ove \O$, so
that $\GM u$ is contained in $A$.
Independently of the choice of $\phi^t$, we can already see that
assumption (b) of Lemma \rf{s3.1.0} is satisfied if (cf.\ (\rf{4.6.1}))
$$
\osc{g} \cdot \sup |\nabla u| \le \alpha \ ,
\eqno (\frm{4.17.1})
$$
while assumption (b') is trivially satisfied, and
$\phi$ has vanishing normal component on $\bd\O\times\R$
(because of (\rf{4.16})). Thus it remains to choose
$\phi^t$ so that (a) and (a') hold, and $\phi$
is approximately regular and divergence-free.

Assumption (a') sets $\phi^t$ equal to $|\nabla u|^2-\beta(u-g)^2$
on the graph of $u$., while requiring that $\phi$ is divergence-free 
in the interior of $A$ yields
$$
\bd_t\phi^t =-\div_x\phi^x=-2\Delta u=-2\beta(u-g)\ .
$$
Integrating in $t$ we obtain that $\phi^t$  is given in $A$ by
$$
\eqalign{
  \phi^t
& =|\nabla u|^2 -\beta (u-g)^2 - 2\beta(u-g)(t-u) \cr
& =|\nabla u|^2 -\beta (t-g)^2 +  \beta (t-u)^2 \ . \cr
}
$$
Therefore assumption (a) of Lemma \rf{s3.1.0}, namely 
$\phi^t\ge |\nabla u|^2-\beta (t-g)^2$, is clearly satisfied in $A$.
Moreover $\phi$ is approximately continuous in $\ove A$ (this is trivial 
for the vectorfield $(\phi^x,0)$, which is continuous, and follows from 
Remark \rf{r1} for $(0,\phi^t)$).

Moreover, $\phi$ is divergence-free in the complement of $A$
if we impose that $\bd_t\phi^t=0$, that is, $\phi^t$ depends only on $x$, 
while the compatibility condition (\rf{2.2.2}), which is required in order
to have that $\phi$ is divergence-free on the entire $\O\times\R$
(cf.\ Remark \rf{s2.9}(b)), yields
$$
\phi^t=\cases{
-\beta\big( {u+M\over 2}-g \big)^2 + \beta\big( {u-M\over 2} \big)^2
     & for $t> {u+M\over 2}$, \cr
-\beta\big( {u+m\over 2}-g \big)^2 + \beta\big( {u-m\over 2} \big)^2
     & for $t< {u+m\over 2}$. \cr
}
$$
Finally one can easily check that $\phi$ is approximately regular
also outside $A$ (see Remark \rf{r1}) and satisfies condition (a) of Lemma
\rf{s3.1.0} as well. 
Therefore $\phi$ is an absolute calibration for $u$, provided that $u$
satisfies (\rf{4.17.1}).

\Rem{\thm{e17.1}}
If inequality (\rf{4.17.1}) is strict, the calibration constructed in
the previous paragraph satisfies the strict inequality in assumption 
(b) of Lemma \rf{s3.1.0}, and by
Remark \rf{s3.5.1} this proves that the solution $u$ of (\rf{4.16}) 
is the unique absolute minimizer of $F$.

\Rem{\thm{e17.2}}
The weaker inequality
$$
\osc{u} \cdot \sup |\nabla u| \le \alpha \ ,
\eqno (\frm{g3})
$$
which follows from (\rf{4.16}) and (\rf{4.17.1}) by the maximum 
principle, is not enough to guarantee the minimality of a solution 
$u$ of (\rf{4.16}), not even the Dirichlet minimality.
Take indeed $n:=1$, $\O:=(-a,a)$, and $g(x):=1$ for $x\ge 0$
$g(x):=-1$ for $x<0$. Then the solution to
(\rf{4.16}) can be computed explicitly:
$$
u(x):= \Big[ 1-{ \cosh(\gamma(1-|x|/a)) \over \cosh\gamma}  \Big] g(x)
\ ,
$$
where $\gamma:=\sqrt\beta\, a$.
Now we fix $\alpha$ so that (\rf{g3}) is satisfied, and precisely
$$
\alpha
:=\osc{u} \cdot \sup |\nabla u| 
 = 2\sqrt\beta\, \tanh\gamma\, \Big[1-{1\over\cosh\gamma}\Big] \ ,
$$
and take the comparison function 
$v(x):=\big[1-1/\cosh\gamma\big] \, g(x)$;
$v$ has the same boundary values as $u$, and a tedious but 
straightforward computation gives
$$
F(v)
=2\sqrt\beta\, \tanh\gamma 
  - {2\sqrt\beta\over \cosh^2\gamma} (\sinh\gamma - \gamma)
\le 2\sqrt\beta\, \tanh\gamma
=F(u) 
\ . 
$$
Therefore $u$ satisfies condition (\rf{g3}) but is not a Dirichlet 
minimizer of $F$.

\Parag{\thm{e18}.~Solution of the Neumann problem for large $\beta$}
The construction in Paragraph \rf{e17} shows that
the solution of the Neumann problem (\rf{4.16}) is an absolute
minimizer of $F$ provided that (\rf{4.17.1}) holds.
However, this condition is far from being necessary.
In particular, for large values of the penalization
parameter $\beta$, the absolute minimizer $u$ of $F$
is close to $g$, and therefore we expect that discontinuities
should not be energetically convenient, at least for
sufficiently regular $g$, and the solution $u$ of
(\rf{4.16}) should be the unique absolute minimizer of~$F$.

We prove this fact by calibration under the assumption that
$\O$ has boundary of class $C^2$, $g$ is of class $W^{2,p}$
for some $p>n$, and $\beta$ is larger than a certain $\beta_0$,  
specified in (\rf{4.26.2}). 
Under these assumptions, $g$ belongs to 
$C^{1,\gamma}(\ove\O)$ for $\gamma:=1-n/p$, 
and $u$ belongs to
$C^{3,\gamma}(\O) \cap C^{1,\delta}(\ove\O) \cap W^{2,p}(\O)$
for every $\delta\in (0,1)$ by the standard regularity theory 
for Neumann problems (see, e.g., [\rf{Tro}], Theorems 3.5, 3.16, 
and 3.17).

Fix a positive constant $\delta$ (to be properly chosen later), 
and take a smooth function $\sigma: \R\to [0,1]$, 
with support included in $[-2\delta,2\delta]$ and identically equal 
to $1$ in $[-\delta,\delta]$, so that $|\dot\sigma|\le 2/\delta$ 
(and then $\| \sigma \|_1\le 4\delta$ and 
$\| \dot\sigma \|_\infty\le 2/\delta$). Set
$$
\phi^x(x,t):=2 \sigma(t-u(x)) \, \nabla u(x) \ .
\eqno (\frm{4.20})
$$
To simplify the notation, in the following we simply write
$\sigma$ and $\nabla u$ instead of $\sigma(t-u(x))$ and $\nabla u(x)$
(this must be kept into account when deriving), so that (\rf{4.20}) 
becomes simply $\phi^x =2\sigma\,\nabla u$.

It follows from (\rf{4.16}) and (\rf{4.20})
that $\phi$ has vanishing normal
component at the boundary of $\O\times\R$, and $\phi^x=2\nabla u$ on the
graph of $u$. Assumption (a') in Lemma \rf{s3.1.0} prescribes
the value of $\phi^t$ on the graph of $u$, and precisely
$$
\phi^t(x,u):=|\nabla u|^2 - \beta(u-g)^2
\quad\hbox{for all $x\in\ove\O$.}
\eqno (\frm{4.21})
$$
We impose now that $\phi$ is divergence-free, which reduces to
$$
\eqalignno{
  \bd_t \phi^t = -\div_x \phi^x
& = -2\sigma \, \Delta u + 2\dot\sigma \, |\nabla u|^2 \cr
& = -2\beta \sigma \, (u-g)  + 2 \dot\sigma \, |\nabla u|^2 \ . 
    & (\frm{4.22}) \cr
}
$$
Identities  (\rf{4.22}) and (\rf{4.21}) together determine
$\phi^t$ everywhere.

Now we want to verify that assumption (a) of Lemma \rf{s3.1.0} holds,
that is,
$$
\textstyle
\phi^t \ge {1\over 4} |\phi^x|^2 - \beta(t-g)^2 \ .
\eqno (\frm{4.22.1})
$$
Since the equality holds by construction on the graph of $u$, 
the full inequality is implied by the following inequalities 
on the derivatives of both sides of (\rf{4.22.1}) with respect 
to $t$:
$$
\cases{
\bd_t\phi^t \ge {1\over 2} \phi^x  \bd_t\phi^x -2\beta(t-g)
   & for $t>u$, \cr
\bd_t\phi^t \le {1\over 2} \phi^x  \bd_t\phi^x -2\beta(t-g) \vspazio
   & for $t<u$. \cr
}
\eqno (\frm{4.23})
$$
Let us consider the first inequality: by (\rf{4.20}) and (\rf{4.22}) it
becomes
$$
-2\beta\sigma\, (u-g)  + 2\dot\sigma\,|\nabla u|^2
\ge 2\sigma \dot\sigma \,|\nabla u|^2 - 2\beta (t-g) \ ,
$$
which is equivalent to
$$
\beta [(t-g) -\sigma(u-g)] \ge \dot\sigma (\sigma-1) |\nabla u|^2 \ .
\eqno (\frm{4.24})
$$
When $u<t\le u+ \delta$ we have $\sigma=1$, and then 
(\rf{4.24}) becomes $t-u\ge 0$, which is obviously true.
When $t>u+ \delta$, we have
$(t-g) -\sigma (u-g) \ge \delta- \|u-g\|_\infty$ and 
$|\dot\sigma (\sigma-1)| \le 2/\delta$,
and then (\rf{4.24}) is implied by
$\beta(\delta- \|u-g\|_\infty)\ge {2\over \delta} \|\nabla u\|_\infty^2$.
This inequality can be rewritten as
$$
\delta^2- \delta \, \|u-g\|_\infty -{2\over \beta}\, 
\|\nabla u \|_\infty^2 \ge 0 \ ,
$$
and is satisfied for
$$
\delta \ge \|u-g\|_\infty + \sqrt {2\over\beta} \, \| \nabla u\|_\infty \ .
\eqno (\frm{4.25})
$$
One checks in the same way that (\rf{4.25}) implies the second
inequality in (\rf{4.23}) too. In other words, assumption (a) of Lemma
\rf{s3.1.0} holds if (\rf{4.25}) holds.

Assumption (b') of Lemma \rf{s3.1.0} is trivially satisfied
because $Su$ is empty, while (\rf{4.20}) and the estimate
$\|\sigma\|_1 \le 4\delta$ imply that assumption (b) of Lemma
\rf{s3.1.0} is satisfied with strict inequality if
$8\delta\| \nabla u\|_\infty < \alpha$, that is, 
$$
\delta < {\alpha\over 8\| \nabla u\|_\infty} \ .
\eqno (\frm{4.26})
$$
Finally, we can find $\delta$ that satisfies both 
(\rf{4.25}) and (\rf{4.26}) if
$$
 \|\nabla u\|_\infty 
\big( \sqrt\beta\, \|u-g\|_\infty + \sqrt 2 \|\nabla u\|_\infty \big)
< {\alpha\over 8} \sqrt\beta\ ,
\eqno (\frm{4.26.0})
$$
and by Theorem \rf{s3.2.0} and Remark \rf{s3.5.1} we conclude that, 
if (\rf{4.26.0}) is satisfied, then $u$ is the {\it unique} absolute 
minimizer of $F$.

Thus it remains to show that (\rf{4.26.0}) holds for $\beta$ large enough.
Note that $u$, being a solution of (\rf{4.16}), depends on $\beta$, 
and there exist positive constants $K$ and $\bar \beta$  
(depending on $\O$, but not on $g$ and $\beta$) such that 
for every $\beta\ge\bar \beta$ there holds
$$
\sqrt\beta\, \| u-g \|_\infty  + \| \nabla u \|_\infty
\le K \| \nabla g \|_{W^{1,p}}
\ .
\eqno (\frm{4.26.1})
$$
This estimate can be derived, for instance, from the interpolation
inequality (3.1.59) of Theorem 3.1.22 in [\rf{Lu}] (one has to replace 
$\lambda$, $u$, ${\cal A}$, and ${\cal B}$ with 
$\beta$, $u-g$, $\Delta$, and $\bd_\nu$ respectively, 
and recall that $\Delta u=\beta(u-g)$). 

Estimate (\rf{4.26.1}) shows that (\rf{4.26.0}) holds
for
$$
\beta> \beta_0:= \max\big\{
\bar\beta, \, 2^7 \alpha^{-2}K^4 \| \nabla g\|_{W^{1,p}}^4
\big\}
\ .
\eqno (\frm{4.26.2})
$$

\Parag{\thm{e19}.~Characteristic functions of regular sets}
If $g:=1_E$ is the characteristic function of a sufficiently regular
compact subset $E$ of $\O$, then it is natural to conjecture that for 
large values of $\beta$ the minimizer of $F$ is $g$ itself.
We prove this statement by calibration under the
assumption that the boundary of $E$ is of class $C^{1,1}$
(cf.\ also Remark \rf{e22} below) and $\beta>\beta_0$, with
$\beta_0$ defined by (\rf{4.33}). Under these assumptions we
also prove the uniqueness of the minimizer.

As in the previous paragraph, we first construct $\phi^x$.
To this end, we take a Lipschitz vectorfield $\psi$ on $\ove\O$
which agrees on $\bd E$ with the inner normal of $\bd E$, is supported
on a neighbourhood of $\bd E$ which is relatively compact in $\O$, and
satisfies $|\psi|\le 1$ everywhere. For instance, we can use that
$\bd E$ is locally a graph, which yields a trivial extension
of the normal vectorfield on a small neighbourhood of each point,
and then use a partition of unity to paste together these
different extensions. Now we set
$$
\phi^x(x,t):=\sigma(t) \, \psi(x)
\quad\hbox{for all $x\in\ove\O$, $t\in\R$,}
\eqno (\frm{4.27})
$$
where $\sigma: \R\to [0,2\alpha]$ is a function of class $C^1$,
supported in $[0,1]$, with integral equal to $\alpha$, and such that
$|\dot\sigma(t)|\le 16 \alpha$ for $t\in[0,1]$, $\sigma(t):=t^2$ for
$t\in[0,1/8]$, $\sigma(t):=(1-t)^2$ for $t\in[7/8,1]$.

We see that, independently of the choice of $\phi^t$,
the vectorfield $\phi$ has vanishing normal component
at the boundary of $\O$, and satisfies assumptions
(b) and (b') of Lemma \rf{s3.1.0}.
Since $\phi^x$ vanishes for $t=0$ and for $t=1$, and therefore on
the graph of $g$, requiring that $\phi$ satisfies assumption (a') yields
$$
\phi^t(x,g(x)):=0
\quad\hbox{for $\L^n$-a.e.\ $x\in\O$,}
\eqno (\frm{4.28})
$$
while requiring that $\phi$ is divergence-free yields
(cf.\ (\rf{4.27}))
$$
\bd_t\phi^t=-\div_x \phi^x=-\sigma \, \div_x \psi \ .
\eqno (\frm{4.29})
$$
Conditions (\rf{4.28}) and (\rf{4.29}) together determine $\phi^t$.

Note that $\phi$ is approximately regular: this is trivial 
for the vectorfield $(\phi^x,0)$, which is continuous, and follows from 
Remark \rf{r1} for $(0,\phi^t)$ (even though $\phi^t$ is discontinuous 
on $\bd E\times\R$ and where $\div_x\psi$ is).

To show that $\phi$ is an absolute calibration
it remains thus to verify assumption (a) of Lemma \rf{s3.1.0}, namely
$$
\textstyle
\phi^t
\ge {1\over 4} |\phi^x|^2 - \beta (t-g)^2 
={1\over 4}  \sigma^2 |\psi|^2 - \beta (t-g)^2 \ .
\eqno (\frm{4.30})
$$
Since the equality holds by construction on the graph of $g$ 
(cf.\ (\rf{4.28})), it is enough that $\bd_t\phi^t$ satisfies the 
inequality
$$
\textstyle
\bd_t\phi^t:=-\sigma \, \div_x \psi
>  {1\over 2} \sigma  \dot\sigma  |\psi|^2 -2\beta (t-g)
\quad\hbox{for $t> g(x)$,}
\eqno (\frm{4.31})
$$
and the opposite inequality for $t<g(x)$.
Inequality (\rf{4.31}) is clearly satisfied for $t>1$, since
$\sigma(t)=0$. If $g(x)=0$ and $0<t\le 1$, (\rf{4.31}) is implied
by
$$
\textstyle
-\sigma\| \div_x \psi\|_\infty > 
{1\over 2}\sigma|\dot\sigma| - 2\beta t\ .
\eqno (\frm{4.32})
$$
In turn, (\rf{4.32}) reduces for $0<t<1/8$ to
$$
-t^2\| \div_x \psi\|_\infty >  t^3 - 2 \beta t \ ,
$$
which is satisfied for
$\beta> {1\over 16} \| \div_x \psi\|_\infty + {1\over 128}$,
while, for $1/8\le t\le 1$, (\rf{4.32}) follows from
$$
\textstyle
-2\alpha \| \div_x \psi\|_\infty
>  16 \alpha^2  -{1\over 4} \beta\ ,
$$
which is satisfied for 
$\beta> 8\alpha\| \div_x \psi\|_\infty + 64\alpha^2$. 
Therefore (\rf{4.31}) holds for
$$
\textstyle
\beta> \beta_0 := \max\Big\{
{1\over 16} \|\div_x\psi\|_\infty + {1\over 128} , 
\ 16\alpha  \|\div_x\psi\|_\infty + 64\alpha^2 \Big\}\ .
\eqno (\frm{4.33})
$$
The same condition implies also the opposite inequality for $t<g(x)$.
This concludes the proof that $\phi$ calibrates $g$.

To prove that $g$ is the unique minimizer of $F$, we first
notice that the strict inequality in (\rf{4.31}) implies the strict 
inequality in (\rf{4.30}) for $t > g(x)$, and of course we have
the strict inequality for $t<g(x)$, too.
In other words, the inequality in assumption (a) of Lemma \rf{s3.1.0}
is strict for all $t\ne g(x)$.
Now, if $u$ is another minimizer, $\phi$ must calibrate $u$, too
(cf.\ Remark \rf{s3.5.1}), and in particular it must satisfy assumption 
(a') of Lemma \rf{s3.1.0} for $u$, which means that
the inequality in assumption (a) is an equality for $t=u(x)$. 
Therefore we conclude that $u(x)=g(x)$ for $\L^n$-a.e.\ $x$ in $\O$. 

\Rem{\thm{e22}}
If $g:=1_E$ is the characteristic function of a set 
$E$ relatively compact in $\O$ and $u:=g$ minimizes $F$, 
then the set $E$ minimizes in particular 
${\Cal F}(A):=F(1_A)=\alpha\H^{n-1}(\bd_* A)+ \beta |A\triangle E|$ 
among all sets $A$ with finite perimeter in $\O$.
Hence the regularity theory for minimal perimeters yields that that, 
in dimension $n\le 7$, $E$ must be of class $C^{1,\gamma}$ for every 
$\gamma<1$, while in dimension two it must be of class $C^{1,1}$
(see, e.g., [\rf{Amb-TGM}], Theorem 4.7.4).
Thus the regularity on $g$ required in the previous paragraph 
is optimal in dimension two, and close to optimal for $3\le n\le 7$.

\Rem{\thm{e22.1}}
If $g$ is the characteristic function of a set $E$ of class $C^{1,1}$
which is not relatively compact in $\O$, then the result of Paragraph 
\rf{e19} can be generalized as follows: $u:=g$ is a minimizer of $F$ 
for large values of $\beta$ provided that $\bd E$ is orthogonal
to $\bd\O$ (which is assumed to be sufficiently smooth). 
An absolute calibration can be constructed as in
the previous paragraph, one has only to choose $\psi$ so that
it is tangent to the boundary of $\O$.

Notice that this orthogonality requirement is necessary:
indeed it is easily proved that, given a minimizer of the functional
$\alpha\H^{n-1}(\O\cap \bd_* A)+ \beta|A\triangle E|$ among all finite 
perimeter sets $A$ in $\O$, its boundary is orthogonal to $\bd\O$.

\medskip
We conclude this section with some remarks on the gradient flow 
associated with the (homogeneous) Mumford-Shah functional.

\Parag{\thm{e23}.~Gradient flow for the Mumford-Shah functional}
A gradient flow for $F_0$ with respect to the $L^2$-metric 
can be defined in a variational way by time discretization, 
following the {\it minimizing movements} approach
developed in [\rf{DG-MM}], [\rf{Amb-MM}], [\rf{Gob2}].
Given an initial datum $u_0\in L^2(\O)$ and a discretization step
$\delta>0$, we set $u_{\delta,0}:=u_0$ and define inductively
$u_{\delta,j}$ for $j=1,2,\ldots$ as {\it any\/} minimizer of
$$
F_0(u)+ {1\over \delta}\int_\O (u-u_{\delta,j-1})^2\, dx
\eqno (\frm{4.35})
$$
among all functions $u$ in $SBV(\O)$ -- with or without prescribed
boundary values, according to the boundary condition (Dirichlet or
Neumann) imposed on the flow. Then we define $u_\delta:
\O\times[0,+\infty)\to\R$ by $u_\delta(x,t):=u_{\delta,j}$ for
$t:=j\delta$, and by linear interpolation for $t\in\big(j\delta,
(j+1)\delta\big)$,
and call {\it gradient flow with initial datum\/} $u(x,0)=u_0(x)$ any 
possible limit of $u_\delta$ as $\delta\to 0$ along any sequence.
Note that the flow may be not unique, as even 
$u_\delta$ is not uniquely defined.

\Rem{\thm{e24}}
If $u_0$ belongs to $W^{1,2}(\O)$, and the minimization of (\rf{4.35})
is restricted a priori to the functions $u$ in $W^{1,2}(\O)$, 
then $F_0(u)$ is just the usual Dirichlet integral,
and it can be proved (cf.\ [\rf{Amb-MM}], Example 2.1)
that the previous construction yields a unique flow
which agrees with the solution of the heat equation 
$$
\bd_t u=\Delta u
\quad\hbox{on $\O\times (0,+\infty)$}
$$
with initial datum $u(x,0)=u_0(x)$ and boundary conditions -- Neumann
or Dirichlet -- according to the boundary conditions imposed in the
minimization of (\rf{4.35}).

\medskip
The previous remark and the result of Paragraph \rf{e18} suggest
that for a smooth initial datum $u_0$, the gradient flow
associated with $F_0$ is just the solution of the heat equation.
To prove this, however, we need some additional information  
on the minima of $F$ for large $\beta$.

\Parag{\thm{e24.1}.~Improved estimates on the solution of (\rf{4.16})}
Under the regularity assumptions on $\O$ and $g$ given in Paragraph 
\rf{e18}, if $\Delta g\in L^\infty(\O)$ and $\bd_\nu g=0$ on $\bd\O$, 
then the solution $u$ to the Neumann  problem {\rm(\rf{4.16})} satisfies
$$
\|\Delta u\|_\infty \le \|\Delta g\|_\infty \ ,
\eqno (\frm{4.26.5})
$$
and an improved version of estimate (\rf{4.26.1}): 
$$
\beta\|u-g\|_\infty + \| \nabla u\|_\infty \le K \|\Delta g\|_\infty \ ,
\eqno (\frm{4.26.6})
$$
with $K$ depending on $\O$, but not on $g$ and $\beta$.
In particular condition (\rf{4.26.0}) of Paragraph \rf{e18} 
holds for 
$$
\beta>\beta_0:=\max\big\{
1, \, 2^7 \alpha^{-2} K^4 \|\Delta g\|_\infty^4\big\} \ ,
\eqno (\frm{4.26.61})
$$
and in that case $u$ is the unique absolute minimizer of $F$.

To prove (\rf{4.26.5}) and (\rf{4.26.6}), we first notice that 
the function $v:=g+\eps$ is a super-solution of (\rf{4.16}) as long as
$\eps\ge\beta^{-1}\|\Delta g\|_\infty$, in the sense that
$$
\cases{
\Delta v\le \beta(v-g) & on $\O$, \cr
\bd_\nu v\le 0 \vspazio & on $\bd\O$. \cr
}
$$
Thus $u$ is (a.e.) smaller than $g+\eps$ on $\O$. 
Similarly, $g-\eps$ is a sub-solution, and then
$$
\|u-g\|_\infty \le {1\over\beta} \|\Delta g\|_\infty \ ,
\eqno (\frm{4.26.7})
$$
which, in view of (\rf{4.16}), implies (\rf{4.26.5}).

Now, $u$ solves the equation $\Delta u=f$ with Neumann boundary 
conditions and $f:=\beta(u-g)$, and well-known estimates (cf. [\rf{Tro}], 
Theorem 3.16) give $\| \nabla u\|_\infty \le K \|f\|_\infty$  
for a suitable constant $K$ depending on $\O$, but not on $f$. 
Together with (\rf{4.26.7}), this implies (\rf{4.26.6}).

\Parag{\thm{e24.3}.~Solution of the Dirichlet problem for large $\beta$}
Assume that $\O$ and $g$ satisfy the regularity assumptions of 
Paragraph \rf{e18}, $\Delta g$ belongs to $L^\infty(\O)$,
and $g_0$ is a function in $W^{2,p}(\O)$, and 
consider the solution $u$ to the Dirichlet problem
$$
\cases{
\Delta u=\beta(u-g) & on $\O$, \cr
u=g_0 \vspazio & on $\bd\O$. \cr
}
\eqno (\frm{4.26.8})
$$
Then $u$ belongs to 
$C^{3,\gamma}(\O)\cap C^{1,\delta}(\ove\O)\cap W^{2,p}(\O)$ 
for every $\delta\in (0,1)$ (see, e.g., [\rf{Tro}], Theorems 3.5, 
3.16, and 3.17).
We claim that if $g_0=g$ on $\bd\O$ and $\beta$ is sufficiently 
large, then $u$ is the unique Dirichlet  minimizer of $F$ with 
boundary value $g_0$. 

This claim can be proved by the same calibration constructed in 
Paragraph \rf{e18}, provided that estimate (\rf{4.26.1}) is suitably
replaced. To this end, we notice that $v:=g+\eps$ is a 
super-solution of (\rf{4.26.8}) 
for $\eps\ge\beta^{-1}\|\Delta g\|_\infty$, in the sense that
$$
\cases{
\Delta v\le \beta (v-g) & on $\O$, \cr
v \ge g_0 \vspazio & on $\bd\O$ \cr
}
$$
(we use here that $g=g_0$ on $\bd\O$), and, similarly, $g-\eps$ is 
a sub-solution. As in the previous paragraph, we deduce that 
$\beta\|u-g\|\le \| \Delta g\|_\infty$, and hence (cf.\ (\rf{4.26.5}))
$$
\| \Delta u\|_\infty \le \| \Delta g\|_\infty \ .
$$
Let us consider now the function $w:=u-g_0$;
since it solves $\Delta w=f-\Delta g_0$ with 
Dirichlet boundary conditions $w=0$ on $\bd\O$, and $f:=\beta(u-g)$, 
well-known estimates for solutions of Dirichlet problems
(see, e.g., [\rf{Tro}], Theorem 3.16) give 
$\| \nabla w\|_\infty \le K (\|f\|_\infty+\| \nabla g_0\|_{C^{0,\gamma}})$.
Together with the estimate on $\|u-g\|_\infty$, this implies
(cf. (\rf{4.26.6}))
$$
\beta \|u-g\|_\infty +\| \nabla u\|_\infty 
\le K \big( \|\Delta g\|_\infty + \| \nabla g_0\|_{C^{0,\gamma}} \big) \ ,
$$
with a possibly different $K$. 
In particular, condition (\rf{4.26.0}) of Paragraph \rf{e18} 
is satisfied for 
$$
\beta>\beta_0:=\max\big\{
1, \, 2^7 \alpha^{-2}K^4 (\|\Delta g\|_\infty
+ \|\nabla g_0\|_{C^{0,\gamma}})^4
\big\} \ ,
\eqno (\frm{4.26.9})
$$
and in that case $u$ is the unique Dirichlet minimizer 
of $F$ with boundary values $u=g_0$.

\Parag{\thm{e25}.~Gradient flow with smooth initial datum}
Assume that $\O$ has boundary of class $C^2$, $u_0\in W^{2,p}(\O)$
for some $p>n$, $\Delta u_0\in L^\infty(\O)$, and $\bd_\nu u_0=0$ 
on $\bd\O$. Then the gradient flow for $F_0$ with Neumann boundary 
conditions and initial datum $u(x,0)=u_0(x)$ constructed in Paragraph 
\rf{e23} is unique, and agrees with the solution of the heat equation.

In virtue of Remark \rf{e24}, this claim is an immediate consequence 
of the following fact: when
$$
\delta < \delta_0:=\Big[
\max\big\{ 1, \, 2^7 \alpha^{-2} K^4 \|\Delta u_0\|_\infty^4 \big\}
\Big]^{-1} \ ,
$$
then, for every integer $j$, {\it every} minimizer of (\rf{4.35}) 
belongs to $W^{1,2}(\O)$. In other words, the solution of the Neumann 
problem (\rf{4.16}) with $\beta:=1/\delta$ and $g:=u_{\delta,j-1}$ 
is the {\it unique} minimizer of (\rf{4.35}).  

To prove this fact, it suffices to verify that the assumptions of 
Paragraph \rf{e24.1} are satisfied for every $j$, and precisely: 
$u_{\delta,j-1}\in W^{2,p}(\O)$, $\Delta u_{\delta,j-1}\in L^\infty(\O)$, 
$\bd_\nu u_{\delta,j-1}=0$ on $\bd\O$, and inequality (\rf{4.26.61}) 
holds with $\beta:=1/\delta$ and $g:=u_{\delta,j-1}$.
The last requirement follows from the choice of $\delta$ and 
the chain of inequalities
$$
\|\Delta u_0\|_\infty:=
\|\Delta u_{\delta,0}\|_\infty \ge \|\Delta u_{\delta,1}\|_\infty \ge
\|\Delta u_{\delta,2}\|_\infty \ge \dots \ ,
$$ 
which are implied by (\rf{4.26.5}).
The $W^{2,p}$ regularity of $u_{\delta,j}$ follows  from the corresponding
regularity of $u_{\delta,j-1}$, as remarked at the beginning of Paragraph 
\rf{e18}.

\Rem{\thm{e25.1}}
The conclusion of the previous paragraph also holds for the 
gradient flow with Dirichlet boundary conditions. More precisely, 
if $\O$ has boundary of class $C^2$ and $u_0$
is of class $W^{2,p}(\O)$, with $\Delta u_0\in L^\infty(\O)$, 
then the gradient flow for $F_0$ with initial datum $u(x,0)=u_0(x)$ on 
$\O$ and Dirichlet boundary condition $u(x,t)=u_0(x)$ on 
$\bd\O\times[0,+\infty)$ constructed  in Paragraph \rf{e23} is unique, 
and agrees with the solution of the heat equation
(with same initial datum and boundary conditions).
The proof is essentially the same as in the Neumann case, 
and relies on the estimates given Paragraph \rf{e24.3}
(in particular one has to apply (\rf{4.26.8}) with
$g=u_{\delta,j-1}$ and $g_0=u_0$).

\Rem{\thm{e25.2}}
If we drop the assumption $\bd_\nu u_0=0$ on $\bd\O$
in Paragraph \rf{e25}, the proof breaks down because we can 
no longer estimate $\| \Delta u_{\delta,j-1}\|_\infty$ by 
$\| \Delta u_0\|_\infty$, but we do not know if the conclusion 
on the gradient flow still holds. 
A similar problem occurs with the statement in Remark \rf{e25.1}
if we replace the Dirichlet boundary condition $u(x,t)=u_0(x)$ 
on $\bd\O$ with a different one.  

\Rem{\thm{e26}}
If $u_0:=1_E$ is the characteristic function of a compact
subset of $\O$ with boundary of class $C^{1,1}$, then the gradient
flow for $F_0$ with Neumann (or Dirichlet) boundary conditions 
constructed in Paragraph \rf{e23} is unique, and
is just given by $u(x,t):=u_0(x)$ on $\O\times[0,+\infty)$. 
This follows immediately from Paragraph \rf{e19}.

\Rem{\thm{e27}}
The conclusions of Paragraph \rf{e25} and Remark \rf{e26} support the
following general conjecture: if the initial datum $u_0$ is smooth out 
of a smooth singular set $Su$, which is compact in $\O$, then, at
least for small times, the associated gradient flow
should be unique and should keep the singular set of $u(\cdot,t)$ equal to
$Su_0$, while the function $u$ evolves in $\O\setminus Su_0$ according to
the heat equation (with Neumann conditions on $Su$). This conjecture
has been proved in the one-dimensional case in [\rf{Gob1}]  (with a
slightly different definition of gradient flow for $F_0$). The general
case will be studied in~[\rf{Mori}] using the calibration method.

\Section{\chp{s6}.~Appendix}
In this section we prove some technical lemmas stated in Section \rf{s2}.
We follow the notation of Section \rf{s2}.

\Lemma{\thm{s6.0}}
{\it
Let $\O$ be an open subset of $\R^n$ whose boundary is the graph of 
a Lipschitz function $f:\R^{n-1}\to\R$, and let $\phi$ be a bounded 
vectorfield on $\ove\O$ which has bounded support and satisfies condition 
{\rm(\rf{2.2.1})} with $M:=\bd\O$. 
Then there exists a sequence of vectors $y_j$ such that $y_j\to 0$, 
$\bd\O+y_j\subset\O$ for every $j$, and
$$
\lim_{j\to\infty} \phi(x+y_j) \cdot \nus{\bd\O}(x) =
\phi(x) \cdot \nus{\bd\O}(x) 
\quad\hbox{for $\H^{n-1}$-a.e.\ $x\in\bd\O$.}
\eqno(\frm{6.1})
$$
}

\Pr
Let $S$ be the set of all vectors $y\in\R^n$ such that 
$\bd\O+y\subset\O$.
For every $r>0$, let $S_r:=S\cap B(0,r)$, and consider the double 
integral
$$
\int_{S_r} \Big[ \int_{\bd\O}
\big| (\phi(x+y) - \phi(x)) \cdot \nus{\bd\O}(x) \big|
\, d\H^{n-1}(x) \Big] \, {dy\over r^n}  \ .
\eqno(\frm{6.2})
$$
If we invert the order of integration, condition {\rm(\rf{2.2.1})} 
means that the inner integral (over $S_r$) tends to $0$ as $r\to 0$ 
for $\H^{n-1}$-a.e. $x\in\bd\O$. Then (\rf{6.2}) converges
to $0$ by the dominated convergence theorem
(recall that $\phi$ is bounded and has bounded support). 

As $\bd\O$ is a graph of a Lipschitz function, the set $S$ contains an 
open cone with vertex in $0$. Then the measure of $S_r$ is larger 
than $ar^n$ for some fixed $a>0$, and
therefore we can choose $y_r\in S_r$ so that the value of the inner 
integral in (\rf{6.2}) is smaller than the double integral divided 
by $a$, and then converges to $0$ as $r\to 0$. 

In other words, $\phi(x+y_r)\cdot\nus{\bd\O}(x)$
converge to $\phi(x)\cdot\nus{\bd\O}(x)$ in $L^1(\bd\O,\H^{n-1})$, 
and then it suffices to choose a subsequence $y_j$ which yields
pointwise convergence for $\H^{n-1}$-a.e. $x\in\bd\O$. \qed

\Parag{Proof of Lemma \rf{s2.0.2}}
(Sketch) We divide the proof in several steps.

\smallskip
{\paragrafo Step 1.} Assume that $\phi$ belong to $C^1_c(\R^n,\R^n)$. 
In this case formula (\rf{2.2.3}) is well-known  -- see, e.g.,
[\rf{Giu}], Theorem 2.10, or [\rf{AFP}], formula (3.87) in Theorem 3.87.

\smallskip
{\paragrafo Step 2.} Assume that $\phi$ is an approximately 
regular vectorfield on $\R^n$ with compact support and that 
$\div\phi\in L^\infty(\R^n)$. 
Let $\psi_\eps(x):=\eps^{-n} \psi(x/\eps)$ be a standard radially 
symmetric regularizing kernel of class $C^\infty_c$, and take 
$\phi_\eps:=\phi*\psi_\eps$.  
Thus formula (\rf{2.2.3}) holds for each $\phi_\eps$ by Step 1, 
and it only remains to check that we can pass to the limit as 
$\eps\to 0$. 
The convergence of the first integral in the right-hand side of 
(\rf{2.2.3}) follows from the fact that the functions $\div\phi_\eps$ 
are bounded in $L^\infty$ and converge to $\div\phi$ a.e.\ in $\O$. 
Since $\phi$ is approximately regular,
the maps $\phi_\eps\cdot\nus{M}$ converge to $\phi\cdot\nus{M}$ 
$\H^{n-1}$-a.e.\ on any Lipschitz surface $M$, 
and then also on any rectifiable set $M$. 
In particular this implies the convergence of the second integral
in the right-hand side of (\rf{2.2.3}). The same argument also 
applies to the left-hand side, provided that we use the coarea formula 
(cf.\ [\rf{AFP}], Theorem 3.40, or [\rf{Fed-GMT}], Theorem 4.5.9(13)) 
to re-write that integral as 
$$
\int_\O \phi_\eps \cdot Du 
=\int_\R \Big[ \int_{M_t} \phi \cdot \nus{M_t} \, d\H^{n-1}\Big] \, dt \ ,
$$
where $M_t$ is the measure theoretic boundary in $\O$ of the sublevel 
$\{u< t\}$. 

\smallskip
{\paragrafo Step 3.} If $\phi$ is a compactly supported, 
approximately regular vectorfield on a neighbourhood of 
$\overline\O$ with  $\div\phi$ in $L^\infty$, we
reduce to Step 2 using a suitable cut-off function.

\smallskip
{\paragrafo Step 4.} Assume that $\O$ is the subgraph of a Lipschitz 
function $f:\R^{n-1}\to\R$, and $\phi$ is a compactly supported,
approximately regular vectorfield  on $\ove\O$ with $\div\phi\in
L^\infty(\O)$.  We take a sequence of vectors $y_j$ as in Lemma 
\rf{s6.0}, and set $\phi_j(x):=\phi(x+y_j)$, $u_j(x):=u(x+y_j)$. 
By Step 3, formula (\rf{2.2.3}) holds with $\phi$ and $u$ replaced by 
$\phi_j$ and $u_j$, and it remains to check that we can pass 
to the limit as $j\to+\infty$.  The convergence is immediate for all 
integrals in (\rf{2.2.3}) but the last one. 
In this case, it suffices to notice that the functions 
$\phi_j\cdot\nus{\bd\O}$ are uniformly bounded  and converge to 
$\phi\cdot \nus{\bd\O}$ $\H^{n-1}$-a.e.\ on $\bd\O$ (by the choice 
of the vectors $y_j$), while the traces of $u_j$ on $\bd\O$ 
converge to the trace of $u$ in $L^1(\bd\O,\H^{n-1})$ (because 
the functions $u_j$ converge to $u$ in variation, or, alternatively, 
because the $L^1$-norm of the difference of the traces is controlled, 
up to a constant which does not depend on $j$, by 
$|Du|(\O\setminus(\ove\O-y_j))$, which clearly tends to zero).

\smallskip
{\paragrafo Step 5.} To prove the general case, we use a locally finite
partition of unity consisting of compactly supported smooth functions 
to reduce to Step 4.  \qed

\Parag{Proof of Lemma \rf{s2.7}}
We first prove that $\div\phi=f$ on $\O\setminus S_0$.
Since the problem is local, it is enough to show that $\div\phi=f$
on every bounded open set $U$ with $\ove U\subset \O\setminus S_0$
and such that $U\setminus S_1$ has two connected
components $U^+$ and $U^-$ with Lipschitz boundary. As $\phi$ is
approximately regular on $\ove U{}^\pm$, we can apply
formula (\rf{2.2.3}) with $\O$ replaced by 
$U^\pm$ and $u\in C^\infty_c(U)$: as the integrals on $U\cap S_1$ cancel 
out, we are left with $\int_U \phi\cdot \nabla u\, dx = -\int_U fu\,dx$. 
By the arbitrariness of $u$, we deduce that $\div\phi=f$ on $U$.

We prove now that $\div\phi=f$ on $\O$. 
Since $\H^{n-1}(S_0)=0$, the $(1,1)$-capacity of
$S_0$ is zero (see [\rf{EG}], Section 5.6.3), and therefore there
exists a sequence of functions $\sigma_j$ in $C^\infty(\ove\O)$
such that $0\le \sigma_j\le 1$ in $\O$, and $\sigma_j=0$ 
in a neighbourhood of $S_0$, and $\sigma_j\to 1$ strongly 
in $W^{1,1}(\O)$.

Take now an arbitrary function $u\in C^\infty_c(\O)$. 
Then the functions $\sigma_j u$ belong to $C^\infty_c(\O\setminus S_0)$ 
and, since $\div\phi=f$ on $\O\setminus S_0$, we have 
$\int_\O \phi\cdot \nabla (\sigma_j u) \, dx=
-\int_\O {f \cdot(\sigma_j u)} \, dx$.
Moreover the functions $\sigma_j u$ converge to $u$ strongly in 
$W^{1,1}(\O)$, and therefore 
$\int_\O  \phi\cdot \nabla u \,  dx=- \int_\O fu\,dx$,
which concludes the proof. 
\qed

\Parag{Proof of Lemma \rf{s2.2}}
By a monotone class argument it is enough to prove (\rf{2.4}) for $\phi$ 
of the form $\phi(x,t):=\rho(t)\, \psi(x)$, with $\rho: \R\to\R$ and
$\psi = (\psi^x,\psi^t):  \O \to \R^n \times \R$ of class $C^\infty_c$.
Let $\sigma$ be the primitive of $\rho$ vanishing at $-\infty$. Then
we have
$$
\eqalignno{
\int_{\O\times\R} \phi \cdot D1_u 
&= - \int_{\O\times\R} \div\phi \, 1_u \, dx \cr 
&= - \int_\O \Big[ \int_{-\infty}^u 
  (\rho\,\div_x\psi^x + \dot\rho \psi^t ) \, dt  \Big] \, dx \cr
&= - \int_\O \big[ \sigma(u)\, \div_x\psi^x + \rho(u)  \psi^t
\big] \, dx\ . & (\frm{6.3}) \cr
}
$$
As $u$ belongs to $SBV(\O)$, the chain-rule for $BV$-functions (see, 
e.g., [\rf{AFP}], Theorem 3.96) gives
$$
D[\sigma(u)]=\rho(u) \nabla u \cdot \L^n+ \big[ \sigma(u^+) - \sigma(u^-)
\big] \nu_u \cdot \H^{n-1} \LL Su \ .
$$
Therefore (\rf{6.3}) implies
$$
\eqalign{
   \int_{\O\times\R} \phi \cdot D1_u 
=& \int_\O \big[ \rho(u) \psi^x \cdot \nabla u - \rho(u)\psi^t \big]
    \, dx \cr
 & +\int_{Su} \big[ \sigma(u^+) - \sigma(u^-) \big] \psi^x \cdot \nu_u 
    \, d\H^{n-1}\ , \cr
}
$$
which, together with (\rf{2.3}), gives (\rf{2.4}) in the case
$\phi(x,t):=\rho(t)\,\psi(x)$.
\qed

\Parag{Proof of Lemma \rf{s2.6}}
We set $w:=1_u-1_v$ on $\O\times\R$. Then $w$ belongs to $BV(\O\times\R)$
and $Dw=\nus{\GM u}\cdot \H^n \LL\GM u- \nus{\GM v}\cdot \H^n \LL\GM v$.

Let us consider the inner trace of $w$ on $\bd U$.
First of all we decompose $\bd U$ as the disjoint union
of $(\O\times\R)\cap\bd U$ and $(\bd\O\times\R)\cap\bd U$.
For every $C^\infty$ vectorfield $\psi$ on $\O\times\R$ with compact 
support we apply formula (\rf{2.2.3}) of Lemma \rf{s2.0.2} with $\O$ 
and $\phi$ replaced by $U$ and $\psi$, respectively, and we obtain
$$
\eqalign{
-\int_{U} w \,\div\psi \, dx 
& = \int_{U} \psi\cdot Dw 
  + \int_{\bd U} w\,\psi\cdot \nus{\bd U}  \, d\H^n \cr
& = \int_{\GM u\cap U} \psi \cdot \nus{\GM u} \, d\H^n
  -\int_{\GM v\cap U} \psi \cdot \nus{\GM v} \, d\H^n \cr
& \quad + \int_{\bd U} w\,\psi\cdot \nus{\bd U}  \, d\H^n\ . \cr
}
\eqno (\frm{app101})
$$
On the other hand, by the definition of distributional derivative we 
have also
$$
\eqalign{
- \int_{\O\times\R} w \,\div\psi \, dx &= \int_{\O\times\R} \psi\cdot Dw 
 \cr
&= \int_{\GM u} \psi \cdot \nus{\GM u} \, d\H^n
-\int_{\GM v} \psi \cdot \nus{\GM v} \, d\H^n\ .
}
\eqno (\frm{app102})
$$
Due to the particular structure
of $U$ and the assumption on the complete graphs of $u$ and $v$,
the function $w$ vanishes a.e.\ on $(\O\times\R)\setminus U$. This 
fact, together with
(\rf{app101}) and (\rf{app102}), implies that the inner trace of $w$ on 
$(\O\times\R)\cap\bd U$ satisfies
$$
w\,\nus{\bd U} = \mathop{1{\trait{0}{0}{2}}}\nolimits_{\GM u} \nus{\GM u} 
- \mathop{1{\trait{0}{0}{2}}}\nolimits_{\GM v} \nus{\GM v}
\qquad \H^n\hbox{-a.e.\ on }(\O\times\R)\cap\bd U\ . 
$$
Therefore $w$ belongs
to $L^1((\O\times\R)\cap\bd U,\H^n)$ and
$$
\int_{(\O\times\R)\cap\bd U} w \,\phi \cdot \nus{\bd U} \, d\H^n =
\int_{\GM u\cap\bd U} \phi \cdot \nus{\GM u} \, d\H^n
-\int_{\GM v\cap\bd U} \phi \cdot \nus{\GM v} \, d\H^n\ .
\eqno (\frm{app1})
$$
Now, the trace of $w$ on $\bd\O\times\R$ is the 
difference of the characteristic functions of the traces of $u$ and 
$v$ on $\bd\O$, and therefore it belongs to
$L^1(\bd\O\times\R,\H^n)$ and
vanishes $\H^n$-a.e.\  on $(\bd\O\times\R)\setminus\bd U$. As
$\nu_{\bd U}=(\nu_{\bd\O},0)$ on $(\bd\O\times\R)\cap\bd U$, this 
implies
$$
\eqalign{
  \int_{(\bd\O\times\R)\cap\bd U} w \,\phi & \cdot \nus{\bd U} \, d\H^n 
  = \int_{\bd\O\times\R} w \,\phi^x \cdot \nus{\bd\O} \, d\H^n \cr
& = \int_{\bd\O} \Big[ \int_v^u \phi^x(x,t) dt\Big] \cdot 
    \nus{\bd\O} d\H^{n-1}\ . \cr
}
\eqno (\frm{app2})
$$
Since the inner trace of $w$ on $\bd U$ belongs to $L^1(\bd U,\H^n)$,
we apply formula (\rf{2.2.3}) of Lemma \rf{s2.0.2} with $\O$ and $u$
replaced by $U$ and $w$, respectively, and get
$$
\int_{\GM u\cap U} \phi \cdot \nus{\GM u} \, d\H^n
- \int_{\GM v\cap U} \phi \cdot \nus{\GM v} \, d\H^n
= -\int_{\bd U} w \,\phi \cdot \nus{\bd U} \, d\H^n \ .
\eqno (\frm{6.5})
$$
Identity (\rf{2.7}) follows now from (\rf{app1}), (\rf{app2}), 
and (\rf{6.5}).
\qed

\Section{References}

\def\refe#1#2#3{\item{\bib{#1}.}\hbox{#2}: #3.}

{\frenchspacing 
\parindent 15 pt
\baselineskip 10 pt
\eightrm

\refe{AAC}
{G. Alberti, L. Ambrosio, X. Cabr\'e}
{On a long-standing conjecture of E. De Giorgi:
symmetry in 3d for general nonlinearities and a local minimality property.
{\eightti Acta Appl. Math.\/}, to appear 
(downloadable from {\eighttt http://cvgmt.sns.it/papers})}

\refe{ABD}
{G. Alberti, G. Bouchitt\'e, G. Dal Maso}
{The calibration method for the Mumford-Shah functional.
{\eightti C. R. Acad. Sci. Paris S\'er. I Math.\/} {\eightbf 329} (1999), 
249--254}

\refe{Amb-ARMA}
{L. Ambrosio}
{Existence theory for a new class of variational problems.
{\eightti Arch. Rational Mech. Anal.\/} {\eightbf 111} (1990), 291--322}

\refe{Amb-MM}
{L. Ambrosio}
{Movimenti minimizzanti.
{\eightti Rend. Accad. Naz. Sci. XL Mem. Mat. (5)\/} {\eightbf 19} (1995), 
191--246}

\refe{Amb-TGM}
{L. Ambrosio}
{{\eightti Corso introduttivo alla teoria geometrica 
della misura ed alle superfici minime\/}.
Appunti dei Corsi Tenuti da Docenti della Scuola 
(Notes of courses given by teachers at the School).
Scuola Normale Superiore, Pisa, 1997}

\refe{AFP}
{L. Ambrosio, N. Fusco, D. Pallara}
{{\eightti 
Functions of bounded variation and free discontinuity problems\/}.
Oxford Mathematical Monographs.
Oxford Science Publications, Oxford, 1999}

\refe{Anz}
{G. Anzellotti}
{{\eightti Traces of bounded vector-fields and the divergence theorem}. 
Unpublished preprint, Dipartimento di Matematica, Universit\`a di Trento, 
Trento, 1983}

\refe{BCDM}
{G. Bellettini, A. Coscia, G. Dal Maso} 
{Compactness and lower semicontinuity properties in 
{\eightti SBD}(\char 10). 
{\eightti Math. Z.} {\eightbf 228} (1998), 337--351}

\refe{Bon-Dav}
{A. Bonnet, G. David}
{Cracktip as a global Mumford-Shah minimizer.
Preprint Univ. Paris Sud, Orsay, 2000 
(downloadable from {\eighttt 
http://www.math.u-psud.fr/$\tilde{}$biblio/pub/2000/})}

\refe{BC}
{G. Bouchitt\'e, A. Chambolle}
{Paper in preparation}

\refe{Bra1}
{K.A. Brakke}
{Minimal cones on hypercubes.
{\eightti J. Geom. Anal.\/} {\eightbf 1} (1991), 329--338}

\refe{Bra3}
{K.A. Brakke}
{Soap films and covering spaces.
{\eightti J. Geom. Anal.\/} {\eightbf 5} (1995), 445--514}

\refe{Bra2}
{K.A. Brakke}
{Numerical solution of soap film dual problems.
{\eightti Experiment. Math.\/} {\eightbf 4} (1995), 269--287}

\refe{Cel}
{P. Celada}
{Minimum problems on {\eightti SBV} with irregular boundary datum. 
{\eightti Rend. Sem. Mat. Univ. Padova} {\eightbf 98} (1997), 193--211}

\refe{Cha}
{A. Chambolle}
{Convex representation for lower semicontinuous functionals in 
{\eightti L}\raise 3 pt\hbox{\fiverm 1}.
{\eightti J. Convex Anal.\/}, to appear
(downloadable from 
{\eighttt http://www.ceremade.dauphine.fr/cadrepub.html})}

\refe{3ecm}
{G. Dal Maso}
{The calibration method for free discontinuity problems.
{\eightti Proceedings of the third european congress of mathematics
(Barcelona, 2000)\/}, to appear}

\refe{DM-M-M}
{G. Dal Maso, M.G. Mora, M. Morini}
{Local calibrations for minimizers of the Mumford-Shah 
functional with rectilinear discontinuity sets.
{\eightti J. Math. Pures Appl. (9)\/} {\eightbf 79} (2000), 141--162}

\refe{DMS}
{G. Dal Maso, J.-M. Morel, S. Solimini}
{A variational method in image segmentation: 
existence and approximation results.
{\eightti Acta Math.\/} {\eightbf 168} (1992), 89--151}

\refe{DG-MM}
{E. De Giorgi}
{New problems on minimizing movements.
{\eightti Boundary value problems for partial differential 
equations and applications (dedicated to E. Magenes)},
81--98, edited by J.-L. Lions and C. Baiocchi.
Research Notes in Applied Mathematics, 29. Masson, Paris, 1993}

\refe{DGA}
{E. De Giorgi, L. Ambrosio}
{Un nuovo funzionale del calcolo delle variazioni.
{\eightti Atti Accad. Naz. Lincei Rend. Cl. Sci. Fis. Mat. Natur. (8)\/}
{\eightbf 82} (1988), 199--210}

\refe{DGCL}
{E. De Giorgi, M. Carriero, A. Leaci}
{Existence theorem for a minimum problem with free discontinuity set.
{\eightti Arch. Rational Mech. Anal.\/} {\eightbf 108} (1989), 195--218}

\refe{EG}
{L.C. Evans, R.F. Gariepy}
{{\eightti Measure theory and fine properties of functions\/}.
Studies in Advanced Mathematics. 
CRC Press, Boca Raton, 1992}

\refe{Fed-GMT}
{H. Federer}
{{\eightti Geometric measure theory\/}.
Grundlehren der mathematischen Wissenschaften, 153. 
Springer-Verlag, Berlin-New York, 1969. 
Reprinted in the series Classics in Mathematics. 
Springer-Verlag, Berlin-Heidelberg, 1996}

\refe{Fed-calib}
{H. Federer}
{Real flat chains, cochains and variational problems.
{\eightti Indiana Univ. Math. J.\/} {\eightbf 24} (1974-75), 351--407}


\refe{Giu}
{E. Giusti}
{{\eightti Minimal surfaces and functions of bounded variation\/}.
Monographs in Mathematics, 80.
Birkh\"auser Boston, Boston, 1984}

\refe{Gob1}
{M. Gobbino}
{Gradient flow for the one-dimensional Mumford-Shah functional.
{\eightti Ann. Scuola Norm. Sup. Pisa Cl. Sci. (4)\/} {\eightbf 27} 
(1998), 145--193}

\refe{Gob2}
{M. Gobbino}
{Minimizing movements and evolution problems in Euclidean spaces.
{\eightti Ann. Mat. Pura Appl. (4)\/} {\eightbf 176} (1999), 29--48}

\refe{Lu}
{A. Lunardi}
{{\eightti Analytic semigroups and optimal regularity in parabolic 
problems\/}.
Progress in Nonlinear Differential Equations and their Applications, 16.
Birkh\"auser Verlag, Basel, 1995}

\refe{Mir}
{M. Miranda}
{Superfici cartesiane generalizzate ed insiemi di perimetro
localmente finito sui prodotti cartesiani.
{\eightti Ann. Scuola Norm. Sup. Pisa Cl. Sci. (3)\/} {\eightbf 18} 
(1964), 515--542}

\refe{M-M}
{M.G. Mora, M. Morini}
{Local calibrations for minimizers of the Mumford-Shah functional 
with a regular discontinuity sets. {\eightti Ann. Inst. H. Poincar\'e, 
Anal. Non Lin\'eaire\/}, to appear
(downloadable from {\eighttt http://www.sissa.it/fa/publications})}

\refe{Mo-Sol}
{J.-M. Morel, S. Solimini}
{{\eightti Variational methods in image segmentation\/}.
Progress in Nonlinear Differential Equations and their Applications, 14.
Birkh\"auser Boston, Boston, 1995}

\refe{Mor}
{F. Morgan}
{Calibrations and new singularities in area-minimizing surfaces: a survey.
{\eightti Variational methods (Paris, 1988)\/}, 329--342, 
edited by H. Berestycki et al.
Progr. Nonlinear Differential Equations Appl., 4.
Birkh\"auser Boston, Boston, 1990}

\refe{Mor2}
{F. Morgan}
{Area-minimizing currents bounded by higher multiples of curves.
{\eightti Rend. Circ. Mat. Palermo (2)\/} {\eightbf 33} (1984), 37--46}

\refe{ML1}
{F. Morgan , G. Lawlor}
{Paired calibrations applied to soap films, immiscible fluids,
and surfaces or networks minimizing other norms.
{\eightti Pacific J. Math.\/} {\eightbf 166} (1994), 55--83}

\refe{Mori}
{M. Morini}
{Global calibrations for the non-homogeneous Mumford-Shah functional.
Paper in preparation}

\refe{MS1}
{D. Mumford, J. Shah}
{Boundary detection by minimizing functionals.
{\eightti IEEE Conference on Computer Vision an Pattern Recognition\/},
San Francisco, 1985}

\refe{MS2}
{D. Mumford, J. Shah}
{Optimal approximation by piecewise smooth functions
and associated variational problems.
{\eightti Comm. Pure Appl. Math.\/} {\eightbf 42} (1989), 577--685}


\refe{Taylor}
{J.E. Taylor}
{The structure of singularities in soap-bubble-like and
soap-film-like minimal surfaces.
{\eightti Ann. of Math. (2)\/} {\eightbf 103} (1976), 489--539}

\refe{Tro}
{G.M. Troianiello}
{{\eightti Elliptic differential equations and obstacle problems\/}. 
The University Series in Mathematics. 
Plenum Press, New York, 1987}

\refe{White}
{B. White}
{The least area bounded by multiples of a curve.
{\eightti Proc. Amer. Math. Soc.\/} {\eightbf 90} (1984), 230--232}

}

{\parindent 0 pt
\baselineskip 10 pt

\vskip 1 cm
\line{%
\vbox{\hsize 5 cm\eightrm
{\autore Giovanni Alberti} \par
Dipartimento di Matematica\par
Universit\`a di Pisa \par
via Buonarroti 2, 56127 Pisa\par
ITALY \par
{\eighttt e-mail: alberti@dm.unipi.it}
}
\vbox{\hsize 5 cm\eightrm
{\autore Guy Bouchitt\'e} \par
UFR des Sciences et Techniques \par
Universit\'e de Toulon et du Var \par
BP 132, 83957 La Garde Cedex \par 
FRANCE \par
{\eighttt e-mail: bouchitte@univ-tln.fr}}
\hfill
}

\bigskip
\vbox{\hsize 5 cm\eightrm
{\autore Gianni Dal Maso} \par
S.I.S.S.A. \par
via Beirut 4, 34014 Trieste\par
ITALY \par
{\eighttt e-mail: dalmaso@sissa.it}}

}

\bye